\documentclass[12pt,reqno,oneside]{amsart}
\usepackage{amsmath,color,bm}
\usepackage{amssymb}
\usepackage{amsthm}
\usepackage{graphics}
\usepackage{eucal}
\usepackage{mathrsfs}
\newtheorem{theorem}{Theorem}[section]

\newtheorem{proposition}{Proposition}[section]

\newtheorem{corollary}{Corollary}[section]

\numberwithin{equation}{section}

\newcommand{\ep}{\varepsilon}

\newcommand{\be}{\mathbf{e}}
\newcommand{\eo}{\mathbf{e}_{0}}

\newcommand{\ea}{\mathbf{e}_{1}}
\newcommand{\ei}{\mathbf{e}_{i}}

\newcommand{\en}{\mathbf{e}_{n}}

\newcommand{\eej}{\mathbf{e}_{J}}

\newcommand{\esj}{\mathbf{e}_{*j}}

\newcommand{\esa}{\mathbf{e}_{*1}}
\newcommand{\ess}{\mathbf{e}_{*\xd}}

\newcommand{\diag}{\operatorname{diag}}

\newcommand{\Span}{\operatorname{span}}
\newcommand{\GL}{\operatorname{GL}}
\newcommand{\Mat}{\operatorname{Mat}}
\newcommand{\bbR}{\mathbb{R}}

\newcommand{\bbZ}{\mathbb{Z}}

\newcommand{\bq}{\mathbf{q}}
\newcommand{\bx}{\mathbf{x}}

\newcommand{\La}{\Lambda}

\newcommand{\Ga}{\Gamma}
\newcommand{\bv}{\mathbf{v}}
\newcommand{\bw}{\mathbf{w}}
\newcommand{\f}{\mathbf{f}}

\newcommand{\ba}{\mathbf{a}}
\newcommand{\Id}{\operatorname{Id}}

\newcommand{\bc}{\mathbf{c}}
\newcommand{\pb}{\pi_{\bullet}}
\newcommand{\oj}{\omega_{j}}

\newcommand{\al}{\alpha}

\newcommand{\cW}{\mathcal{W}}

\newcommand{\xs}{s}
\newcommand{\xd}{d}

\newcommand{\R}{{\mathbb R}}
\newcommand{\Z}{{\mathbb Z}}

\newif\ifdraft

\newif\ifcolorcomments

\colorcommentstrue

\newcommand{\allowcomments}[4]{
\newcommand{#1}[1]{\ifdraft{\ifcolorcomments{\textcolor{#4}{##1 --#3}}\else{\textsl{##1 \ --#3}}\fi}\else{}\fi}
}

\allowcomments{\comvictor}{VB}{Victor}{blue}
\allowcomments{\comanish}{AG}{Anish}{blue}
\allowcomments{\comsanju}{SV}{Sanju}{blue}

\begin{document}
\title[Inhomogeneous dual Diophantine approximation]{Inhomogeneous dual Diophantine approximation on affine subspaces}

\begin{abstract}
We prove the convergence and divergence cases of an inhomogeneous Khintchine-Groshev  type theorem for dual approximation restricted to affine subspaces in $\bbR ^n$. The divergence results are proved in the more general context of Hausdorff measures.
\end{abstract}
\author{Victor Beresnevich}
\author{Arijit Ganguly}
\author{Anish Ghosh}
\author{Sanju Velani}
\thanks{VB and SV are supported by EPSRC Programme grant EP/J018260/1. Ghosh is supported by UGC}

\address{School of Mathematics, Tata Institute of Fundamental Research, Mumbai, 400005, India}
\email{arimath@math.tifr.res.in, ghosh@math.tifr.res.in}
\address{Department of Mathematics,  University of York, Heslington, York YO10 5DD, UK}
\email{victor.beresnevich@york.ac.uk, sanju.velani@york.ac.uk}

\maketitle
\section{Introduction}

Throughout the paper,  $\psi : \bbR_{+} \to \bbR_{+} $ is  a non-increasing function and  $\cW_{n}(\psi)$ is  the set of $\bx \in \bbR^n$ for which there exist infinitely many $\ba \in \mathbb{Z}^n\setminus \{\mathbf{0}\}$ such that
\begin{equation}\label{preKG}
|a_0 +  \bx \cdot \ba| < \psi(\|\ba\|^n)
\end{equation}
\noindent for some $a_0 \in \bbZ$. Here and throughout,   $\|~\|$ denotes the supremum norm of a vector and the dot stands for the standard inner product of vectors. For obvious reasons, the set $\cW_{n}(\psi)$ is often referred to as the (dual) set of ``$\psi $-approximable" vectors in $\R^n$. The fundamental Khintchine-Groshev theorem \cite{Groshev, Khintchine} in the metric theory of Diophantine approximation,  provides  an elegant  characterization of the $n$-dimensional Lebesgue measure of $\cW_{n}(\psi)$ in terms of the  convergence/divergence  properties of a `volume sum' associated with the \emph{approximating} function $\psi$. We reinforce the fact that  $\psi$ will always be assumed to be non-increasing.

\begin{theorem}[Khintchine-Groshev] \label{KG}  Let $\psi$ be an  approximating function.  Then
\begin{equation}
|\cW_{n}(\psi)| = \left\{
\begin{array}{rl}
0 & \text{if }  \ \sum_{k=1}^\infty \psi(k) < \infty\\
\\
\text{ full } & \text{if } \  \sum_{k=1}^\infty  \psi(k) = \infty.
\end{array} \right.
\end{equation}
\end{theorem}

\noindent We will use $|~\,\,|$ to denote the absolute value of a real number as well as the Lebesgue measure of a measurable subset $X$  of $\bbR^n$; the context will make the use clear.  

\vspace*{1ex}

\noindent {\em Remarks}\/:
\begin{enumerate}
\item By Dirichlet's theorem, $\cW_{n}(\psi) =  \bbR^n  $  when $\psi(k)=k^{-1}$.
\vspace*{2ex}

\item  A  point $\bx \in \bbR^n$  is called \emph{very well approximable} (abbr. VWA) if there exists $\ep >0 $ such that   $\bx   \in \cW_{n}(\psi_{\ep})
$ where $$ \psi_{\ep} : \bbR_{+} \to \bbR_{+} : k \to  \psi_{\ep}(k):=k^{-(1+ \ep)}  \, . $$ Thus, the essence of the definition of VWA points is that for these points the `Dirichlet exponent'  can be improved beyond the trivial.   Note that in view of  Theorem \ref{KG}, we have that $|\cW_{n}(\psi_{\ep})|= 0$ for any $\ep > 0$.  In other words, almost every point $\bx \in \bbR^n$  is not VWA.

\vspace*{2ex}

\item  The more general Hausdorff measure version of the Khintchine-Groshev theorem has been established in \cite{det}.   For a general background to the classical theory of metric Diophantine approximation,  we refer the reader to the survey type articles \cite{Beresnevich-Bernik-Dodson-Velani-Roth, VicFel}.
\end{enumerate}

We now consider the setting of inhomogeneous Diophantine approximation. Let  $\theta : \mathbb{R}^n \to \mathbb{R}$ be a  function and, given $\psi$,  we let   $\cW^{\theta}_{n}(\psi)$ be the set of $\bx \in \bbR^n$ for which there exist infinitely many $\ba \in \mathbb{Z}^n\setminus \{\mathbf{0}\}$ such that
\begin{equation}\label{preKGinhom}
|a_0 +  \bx \cdot \ba + \theta(\bx)| < \psi(\|\ba\|^n)
\end{equation}
\noindent for some $a_0 \in \bbZ$.   The set $\cW^{\theta}_{n}(\psi)$ is often referred to as the (dual) set of ``$(\psi, \theta)$-inhomogeneously approximable" vectors in $\R^n$.  The following inhomogeneous version of Theorem~\ref{KG} is established in \cite{BaBeVe}. We denote by $C^n$ the set of $n$-times continuously differentiable functions.

\begin{theorem}\label{KGinhom} Let $\psi$ be an  approximating function and $\theta : \mathbb{R}^n \to \mathbb{R}$  be a function  such  that $\theta  \in C^2$.  Then
\begin{equation}
|\cW^{\theta}_{n}(\psi)| = \left\{
\begin{array}{rl}
0 & \text{if } \ \sum_{k=1}^\infty  \psi(k) < \infty\\
\\
\text{ full } & \text{if } \ \sum_{k=1}^\infty  \psi(k) = \infty.
\end{array} \right.
\end{equation}
\end{theorem}

\noindent We remark that the choice of $\theta = \text{constant}$ is the setting of traditional inhomogeneous Diophantine approximation and in that case the above result was well known, see  for example \cite{Cassels}. For the more general Hausdorff measure version of Theorem~\ref{KGinhom} within the traditional setting  see \cite[\S12.1]{Beresnevich-Dickinson-Velani-06:MR2184760, VicFel}.

In this paper, we consider the theory of  Diophantine approximation on manifolds, specifically inhomogeneous approximation on  affine subspaces. The subject of metric Diophantine approximation on manifolds studies the conditions under which a smooth submanifold of $\bbR^n$ inherits Diophantine properties of $\bbR^n$ which are generic for Lebesgue measure. Examples include the resolution of the famous Baker-Sprind\v{z}uk conjecture \cite{Sp} due to Kleinbock and Margulis \cite{KM} using homogeneous dynamics on the space of unimodular lattices. Their result implies that almost every point on a nondegenerate submanifold  $\mathcal{M}$ of $\bbR^n$ is not very well approximable; that is \\
\begin{equation}   \label{KME}
 |\cW_{n}(\psi_{\ep})  \cap   \mathcal{M} \ |_{\mathcal{M}}= 0   \quad  \forall \ \ \ep > 0  \, .
\end{equation}

\noindent Here  and elsewhere $| \ . \ |_{\mathcal{M}}$ denotes the induced Lebesgue measure  on
$\mathcal{M}$. It is worth mentioning that any manifold  $\mathcal{M}$ of $\bbR^n$  satisfying  \eqref{KME}  is called extremal and that  Kleinbock and Margulis proved the stronger ``multiplicative" extremal statement for nondegenerate manifolds.  Essentially, nondegenerate manifolds are smooth manifolds of  $\R^n$ which are sufficiently
curved so as to deviate from any hyperplane, see \cite{Ber1, KM} for a formal definition.

The convergence case of the Khintchine-Groshev theorem was shown to hold for nondegenerate submanifolds of $\R^n$ in \cite{BKM} and independently in \cite{Ber1}. Indeed, in \cite{BKM} Bernik, Kleinbock and Margulis established the stronger `multiplicative' version. The complementary divergence case was subsequently proved in \cite{BBKM} and, as a result, we have the following complete analogue of the Khintchine-Groshev theorem  for nondegenerate manifolds.

\begin{theorem}\label{BBKMthm}  Let $\mathcal{M}$ be a
 nondegenerate submanifold  of $\bbR^n$ and  let
 $\psi$ be an  approximating function. Then
 \begin{equation}  \label{BBKMthmeq}
|\cW_{n}(\psi)  \cap   \mathcal{M} \ |_{\mathcal{M}}
= \left\{
\begin{array}{rl}
0 & \text{if } \ \sum_{k=1}^\infty  \psi(k) < \infty\\
\\
\text{ full } & \text{if } \ \sum_{k=1}^\infty  \psi(k) = \infty.
\end{array} \right.
\end{equation}
\end{theorem}

 \noindent We note that the convergence case of above theorem  implies the extremal statement \eqref{KME} for nondegenerate manifolds.  In this paper, we are concerned with affine subspaces which are the main examples of manifolds which are not nondegenerate.   Theorem \ref{BBKMthm} above is therefore not applicable to them. Nevertheless,  the analogue of the Baker-Sprind\v{z}uk conjecture for affine subspaces was studied by Kleinbock in \cite{Kleinbock-extremal} (see also \cite{Kleinbock-exponent}), and the Khintchine-Groshev theorem in the series of  papers \cite{BBDD, G1, G-mult, G-Monat}. We refer the reader to the recent survey \cite{G-handbook} for further details on this subject.   The key goal of this paper is to  develop the analogous inhomogeneous  theory for affine subspaces.

We now briefly describe the current state of the inhomogeneous theory of  Diophantine approximation on manifolds.  In \cite{BeVe, BeVe2},  the authors discovered a transference principle which allowed them to establish the  inhomogeneous version of the Baker-Sprind\v{z}uk conjecture for nondegenerate manifolds  from the original homogeneous  statement.  Indeed, the inhomogeneous ``multiplicative'' extremal  statement established in \cite{BeVe} implies that for any  nondegenerate  submanifold  $\mathcal{M}$ of $\bbR^n$  and  $\theta = \text{constant},$
\begin{equation}   \label{KMEinhom}
 |\cW^{\theta}_{n}(\psi_{\ep})  \cap   \mathcal{M} \ |_{\mathcal{M}}= 0   \quad  \forall \ \ \ep > 0  \, .
\end{equation} \\
It is worth mentioned that in \cite{GM}, it has been shown that the homogeneous to inhomogeneous transference principle of \cite{BeVe} is flexible enough to be used for arbitrary Diophantine exponents, not just the critical or ``extremal" one.  As demonstrated in \cite{GM},  this naturally extends the scope of potential applications of the original transference principle.   Beyond extremal statements such as \eqref{KMEinhom},  the complete inhomogeneous version of the Khintchine-Groshev theorem, both convergence and divergence cases, for nondegenerate manifolds is established in \cite{BaBeVe}.  In other words,  with mild conditions imposed on the `inhomogeneous'  function $\theta,$ the statement of  Theorem \ref{BBKMthm} is shown to be valid with $\cW_{n}(\psi)$ replaced by $\cW^{\theta}_{n}(\psi)$.    In fact, in  the divergence case, for any $\theta  \in C^2$ the more general Hausdorff measure version is established. As is to be expected,  the convergence case of the Khintchine-Groshev theorem established in \cite{BaBeVe} implies the  inhomogeneous  extremal statement \eqref{KMEinhom} for nondegenerate manifolds.  To the best of our knowledge, unlike  in the homogeneous setting,  an  inhomogeneous theory of  Diophantine approximation on affine subspaces is yet to be developed.  As already alluded to above,  the purpose of this  work is to address this imbalance by establishing  an inhomogeneous version of the Khintchine-Groshev theorem for affine subspaces of $\R^n$.  As a consequence, we obtain the  inhomogeneous extremal statement \eqref{KMEinhom} for affine subspaces.   Indeed, our results go some way towards developing a coherent inhomogeneous theory for degenerate manifolds as outlined in \cite[\S1.4]{BaBeVe}.


%

In the study of Diophantine approximation on affine subspaces, one needs to assume some condition on the slope of the affine subspace in order to ensure that the affine subspace inherits generic Diophantine properties from its ambient Euclidean space. We will now introduce certain Diophantine exponents of matrices which play a key role in this regard.  Indeed, we need these exponents in order to even   state our main convergence theorem.

\subsection{Diophantine exponents of
matrices}\label{sec:exp}

Throughout $\mathscr{H}$ will be an open subset of a $\xd $-dimensional affine subspace of
$\bbR^n$. By making a change of variables if necessary we can assume without loss of generality that $\mathscr{H}$ is
of the
form
\begin{equation}\label{vb5}
\{ (\mathbf{x},\mathbf{x} A' + \mathbf{a}_0): \mathbf{x}\in U\},
\end{equation}
where $\mathbf{a}_0 \in
\mathbb{R}^{n-\xd }$ and $A' \in \Mat_{\xd \times n-\xd }(\bbR)$ and $U$ is an open subset of $\bbR^\xd$.
On setting
 $$A  := \begin{pmatrix}\mathbf{a}_0\\A' \end{pmatrix}  \, , $$ we can rewrite this parametrisation as
\begin{equation}\label{defaffine}
\bx \mapsto (\bx, \tilde{\bx}A)\quad \text{where }  \ \tilde{\bx} :=
(1, \bx).
\end{equation}

\noindent Given a column $\bm\theta\in\R^{\xd+1}$ and a matrix
$A \in \Mat_{\xd+1 \times n-\xd}(\bbR)$, the inhomogeneous Diophantine exponent $\omega(A;\bm\theta)$ of $(A;\bm\theta)$ is defined to be the
supremum of $v > 0$ for which there are infinitely many $\ba' \in \mathbb{Z}^{n-\xd}\setminus \{\mathbf{0}\}$  such that
\begin{equation}\label{defexponent}
\|A \ba' + \ba''+\bm\theta\| < \|\ba'\|^{-v}
\end{equation}
\noindent for some $\ba'' \in \bbZ^{\xd+1}$. In the case $\bm\theta=\bm0$, $\omega(A):=\omega(A;\bm0)$ is the usual (homogeneous) Diophantine approximation exponent of the matrix $A$. It is well known that
$(n-\xd)/(\xd+1) \leq \omega(A) \leq \infty$
for all $A \in \Mat_{\xd+1 \times n-\xd}(\bbR)$ and that $\omega(A)
= (n-\xd)/(\xd+1)$ for Lebesgue almost
every $A$.

We now introduce the higher Diophantine
exponents of $A$ as defined by Kleinbock in
\cite{Kleinbock-exponent}. For $A \in \Mat_{\xd +1 \times
n-\xd }(\bbR)$, we set
\begin{equation}\label{defR}
R_A := \begin{pmatrix}\Id_{\xd +1} & A \end{pmatrix}~\in~ \Mat_{\xd +1 \times
n+1 }(\bbR)
\end{equation}
where $\Id_{\xd+1} $  denotes the $(\xd+1)  \times (\xd+1)$ identity matrix. Let $\eo, \dots, \en$ denote the standard
basis of $\bbR^{n+1}$ and set
\begin{equation}
W_{i \to j} := \Span\{\be_i, \dots, \be_j\},
\end{equation}
where $i\le j$, be the linear subspace of $\R^{n+1}$ spanned by vectors $\be_i, \dots, \be_j$.
Clearly, $W_{0 \to n}=\R^{n+1}$. Now let $\bw \in \bigwedge^{j}(W_{0 \to n})$ represent
a discrete subgroup $\Gamma$ of $\bbZ^{n+1}$, that is $\bw$ is the wedge product of vectors from any given basis of $\Gamma$. Define
the map

 $$\bc : \bigwedge\nolimits^{j}(W_{0 \to n}) \to
(\bigwedge\nolimits^{j-1}(W_{1 \to n}))^{n+1}$$ by setting
\begin{equation}\label{defc}
\bc(\bw)_i := \sum_{\substack{J \subset \{1, \dots, n\}\\ \#J
= j-1}} \langle \be_i \wedge \eej, \bw\rangle \eej
\end{equation} \\
for $0\le i\le n$, and let $\pb$ denote the projection
$\bigwedge(W_{0 \to n}) \to \bigwedge(W_{\xd +1 \to n})$. For
each $j = 1, \dots, n-\xd $, define
\begin{equation}\label{defexponenthigher}
\oj(A) := \sup\left\{ v:\left. \aligned \exists\,
\bw \in \bigwedge\nolimits^{j}(\bbZ^{n+1}) \text{ with arbitrary
large
   } \|\pb(\bw)\|  \\[1ex]
\text{   such that   }\|R_{A} \
\bc(\bw)\|
< \|\pb(\bw)\|^{-\frac {v+1-{j}}{j}}\ \
\endaligned\right.\right\}.
\end{equation}

\noindent It is shown in Lemma 5.3 of \cite{Kleinbock-exponent} that $\omega_1(A) = \omega(A)$.    We shall see in the next section  that the Diophantine exponents $\oj(A)$  play a key role in the  convergence case of the Khintchine-Groshev theorem for affine subspaces.

\subsection{Our Main Theorems}
As before $\mathscr{H}$  will denote an open subset of a $\xd $-dimensional affine subspace of $\R^n$
parametrised as in (\ref{defaffine}).  Then, given  $\psi$  and  $\theta $, the object of study is the set $\cW^{\theta}_{n}(\psi) \cap \mathscr{H}$; i.e.  the set of ``$(\psi, \theta)$-inhomogeneously approximable" vectors on $\mathscr{H}$. \footnote{The reason for considering an open subset of an affine subspace rather than the whole subspace is to allow inhomogeneous functions $\theta$ that may not necessarily be defined on the whole subspace, for example, $\theta(\bx)=\sqrt{1-(x_1^2+\cdots+x_\xd^2).}$} Our first result establishes the convergence case of the inhomogeneous  Khintchine-Groshev theorem for affine subspaces.
\begin{theorem}\label{main convergence}
Let $\mathscr{H}$  be  an open subset of an affine subspace of $\R^n$ of dimension $\xd $ given by \eqref{defaffine}, and suppose that
\begin{equation}\label{diocond}
\oj(A) < n~\text{for every}~j = 1, \dots, n-\xd.
\end{equation}
Let $\psi$ be an  approximating function and $\theta : \mathbb{R}^n \to \mathbb{R}$  be a function  such  that $\theta|_{\mathscr{H}}$ is analytic. Further in the case $\theta|_{\mathscr{H}}$ is a linear function so that
\begin{equation}\label{hatteta}
\hat\theta(\bx):= \theta(\bx, \tilde{\bx}A)=\tilde{\bx}\bm\theta=\theta_0+\theta_1x_1+\dots+\theta_\xd x_\xd
\end{equation}
for some column $\bm\theta=(\theta_0,\dots,\theta_n)^t$, assume that
\begin{equation}\label{diocond+}
\omega(A;\bm\theta)<n\,.
\end{equation}
Then,
\begin{equation}\label{eq:final}
|\cW^{\theta}_{n}(\psi) \cap \mathscr{H}|_{\mathscr{H}}=0
\end{equation}
whenever
\begin{equation}\label{eqn:conv1}
\displaystyle \sum_{k=1} ^{\infty} \psi(k)\, <
\,\infty\,.
\end{equation}
\end{theorem}

\noindent Recall that  $|~|_{\mathscr{H}}$ denotes the  induced Lebesgue measure on $\mathscr{H}$.  Clearly the above convergence theorem implies  the  inhomogeneous extremal statement \eqref{KMEinhom} for any affine subspace satisfying \eqref{diocond} and any analytic $\hat\theta$ that additionally satisfies \eqref{diocond+} in the case it is linear.

\bigskip

\noindent {\em Remarks}\/: 
\begin{enumerate}
\item In the case $\theta|_{{\mathscr{H}}}$ is linear, condition \eqref{diocond+} on the exponent of $(A;\bm\theta)$ is optimal. Indeed, suppose that
\begin{equation}\label{defexponent+}
\|A \ba' + \ba''+\bm\theta\| < \|\ba'\|^{-n}\log^{-3}\|\ba'\|
\end{equation}
holds for infinitely many $\ba'\in\mathbb{Z}^{n-\xd}\setminus \{\mathbf{0}\} $ and some $\ba''\in\Z^{\xd+1}$, but
\begin{equation}\label{defexponent++}
\|A \ba' + \ba''+\bm\theta\| < \|\ba'\|^{-n}\log^{-4}\|\ba'\|
\end{equation}
holds only for finitely many $\ba'\in\mathbb{Z}^{n-\xd}\setminus \{\mathbf{0}\}$ and $\ba''\in\Z^{\xd+1}$. Clearly, in this case $\omega(A;\bm\theta)=n$. The existence of such pairs $(A;\bm\theta)$ can be proved by using the inhomogeneous version of Jarnik's theorem for systems of linear forms \cite[Theorem~19]{Beresnevich-Bernik-Dodson-Velani-Roth}. Assuming $\ba'$ and $\ba''$ satisfy \eqref{defexponent+}, write $\ba'$ as $(a_{\xd+1},\dots,a_n)^t$ and $\ba''$ as $(a_0,\dots,a_\xd)^t$. Then, on identifying $\ba$ with $(a_1,\dots,a_n)^t$ one readily verifies that
\begin{equation}\label{vb1}
\hat{\theta}(\bx) + a_0 + (\bx,
\tilde{\bx}A)\ba=\tilde{\bx}(A \ba' + \ba''+\bm\theta)
\end{equation}
and that $\|\ba\|\ll \|\ba'\|$. Therefore
\begin{align*}
|\hat{\theta}(\bx) + a_0 + (\bx,\tilde{\bx}A)\ba|&\ll \|A \ba' + \ba''+\bm\theta\|\ll \|\ba\|^{-n}\log^{-3}\|\ba\|\,.
\end{align*}
Take $\psi(h)=h^{-1}(\log h)^{-2}$. Then, clearly for every $\bx$ the inequality
\begin{align*}
|\hat{\theta}(\bx) + a_0 + (\bx,\tilde{\bx}A)\ba|<\psi(\|\ba\|^n)
\end{align*}
holds for infinitely many $\ba\in\mathbb{Z}^{n-\xd}\setminus \{\mathbf{0}\}$ and some $a_0\in\Z$. In this case
$$
\cW^{\theta}_{n}(\psi)  \cap   \mathscr{H} =\mathscr{H}
$$
despite \eqref{eqn:conv1}.
An obvious modification of the above argument shows that, given an approximating function $\psi$ satisfying \eqref{eqn:conv1}, in the case $\theta|_{\mathscr{H}}$ is linear, \eqref{eq:final} necessarily implies the existence of $c > 0$ such that
\begin{align*}
\qquad\|A \ba' + \ba''+\bm\theta\|\ge c\ \psi(\|(\ba',\ba'')\|^{n})\qquad\text{for all $(\ba',\ba'')\in \mathbb{Z}^{n-\xd}\setminus \{\mathbf{0}\}\times \Z^{\xd+1}$.}
\end{align*}
How close this is to being a sufficient condition remains an interesting question that is open even in the homogeneous case.\\
\item We note that in case $\bm\theta=\bm0$, since $\omega(A):=\omega(A;\bm0)$ the inhomogeneous Diophantine condition (\ref{diocond+}) does not add extra hypotheses in the homogeneous case.
\end{enumerate}

\bigskip

For the divergence counterpart to Theorem  \ref{main convergence}, we shall prove the following more general statement in terms Hausdorff measures.   Throughout, $\mathcal{H}^\xs  (X) $ will denote the $\xs$-dimensional Hausdorff measure of a subset $X$ of $\R^n$ and $\dim X$ the Hausdorff dimension, where $\xs>0$ is a real number.

\begin{theorem}\label{main divergence}
Let $\mathscr{H}$ an open subset of an affine subspace of $\R^n$ of dimension $\xd $ and let $\xs >\xd-1$.
Let $\psi$ be an  approximating function and $\theta : \mathbb{R}^n \to \mathbb{R}$  be a function  such  that $\theta|_{\mathscr{H}}\in~C^2$.  Suppose that \eqref{diocond} holds and that
\begin{equation}\label{eq:divcond}
 \displaystyle \sum_{k=1}^\infty
 k^{\frac{\xd-\xs}n}\psi (k)^{\xs+1-\xd }=\infty\ .
\end{equation}
Then
\begin{equation}\label{eq:final div}
\mathcal{H}^\xs(\cW^{\theta}_{n}(\psi) \cap \mathscr{H})= \mathcal{H}^\xs(\mathscr{H})\,.
\end{equation}
\end{theorem}

\vspace{2ex}

Given an approximating function $\psi$, the lower order at infinity $\tau_{\psi}$ of $1/\psi$ is defined by
\begin{equation}
\tau_{\psi} := \liminf_{t \to \infty}\frac{-\log\psi(t)}{\log t}
\end{equation}
and indicates the growth of  $1/\psi$ `near' infinity.    Now observe that the divergent sum condition  \eqref{eq:divcond}  is satisfied whenever $$\xs < \xd - 1 + (n + 1)/(n\tau_{\psi} + 1).$$  Therefore, it follows from the definition of Hausdorff dimension that

$$\dim(\cW^{\theta}_{n}(\psi) \cap \mathscr{H})\geq \xs   \quad \text { if }   \quad
 \mathcal{H}^\xs(\cW^{\theta}_{n}(\psi) \cap \mathscr{H}) > 0$$ \\    and that $\mathcal{H}^\xs(\mathscr{H}) > 0$ if $\xs \leq \dim \mathscr{H} = \xd$ and $\mathscr{H}\neq\emptyset$. We therefore obtain the following dimension statement  concerning  the set $  \cW^{\theta}_{n}(\psi) \cap \mathscr{H}  $.  \\

 \begin{corollary}
 Let $\mathscr{H}$ be a non-empty open subset of an  affine subspace of $\R^n$ of dimension $\xd $.
Let $\psi$ be an  approximating function and $\theta : \mathbb{R}^n \to \mathbb{R}$  be a function  such  that $\theta|_{\mathscr{H}}\in~C^2$.  Suppose that \eqref{diocond} holds and  that $1 \leq \tau_{\psi} < \infty$. Then
 \begin{equation}
 \dim(\cW^{\theta}_{n}(\psi) \cap \mathscr{H}) \geq \xd - 1 + \frac{n + 1}{n\tau_{\psi} + 1}.
 \end{equation}
 \end{corollary}

\noindent {\em Remarks}\/:
\begin{enumerate}
\item To the best of our knowledge, the above findings constitute the first known results in the context of inhomogeneous Diophantine approximation on affine subspaces. In fact, Theorem \ref{main divergence} is new even for Lebesgue measure (i.e. when $\xs=\xd$) and in many cases for the homogeneous setting  (i.e. when $ \theta \equiv 0 $). The only previously known cases in the homogeneous setting were:\\[-1ex]
\begin{enumerate}
\item the case of lines passing through the origin which was treated in \cite{BBDD}, and \\[-1ex]
 \item the case of affine hyperplanes ($\xd =n-1$) which was treated in \cite{G-div}. \vspace*{2ex}
\end{enumerate}

\vspace*{2ex}

\item
In the case $\xd =n-1$ condition \eqref{diocond} represents a single inequality imposed on the main Diophantine exponent $\omega(A)$ of $A$.

\vspace*{2ex}

\item In the case $\xs = \xd$\/ the sum within \eqref{eq:divcond} matches the one within \eqref{eqn:conv1}. Thus, Theorem~\ref{main divergence} naturally  complements the statement of Theorem~\ref{main convergence}.

\vspace*{2ex}
\item  In \cite{BaBeVe}, the smoothness condition imposed on the inhomogeneous function $\theta$ is weaker than what we have assumed to establish the convergence statement of Theorem \ref{main convergence}.  In short we  have imposed the  stronger analyticity condition to deal with a technical problem involving $(C, \alpha)$-good functions (see below for the definition). It is plausible that this condition can be relaxed and brought at par with that imposed in \cite{BaBeVe}.

 \vspace*{2ex}

\item  
    The homogeneous results in \cite{BKM} and the inhomogeneous convergence results in \cite{BaBeVe} are proved in the context of more general multivariable approximating functions. This setting includes the case of ``multiplicative" Diophantine approximation.  Both our main theorems should hold for nondegenerate submanifolds of affine subspaces and our convergence theorem should, in addition, be true in the multivariable setting. We plan to return to this extension in a separate work.

\end{enumerate}

\section{The Gradient Division}\label{gradientdivision}

In this section, we prepare the groundwork to prove Theorem~\ref{main convergence}, the ``convergence case". Let $U$ be the same as in \eqref{vb5} and as before define   $\hat{\theta} :U\longrightarrow \mathbb{R}$ by setting
$$\hat\theta(\bx):=\theta ((\bx,\tilde{\bx}A)).$$ Clearly $\hat{\theta}$ is an analytic function since $\theta|_{\mathscr{H}}$ is analytic. For $\mathbf {a} \in
\mathbb{Z}^n\setminus \{\mathbf{0}\}$, we define
\[\mathcal {L}(\mathbf {a}):= \{ \mathbf{x}\in U: |\hat{\theta}(\bx) + a_0 + (\bx,
\tilde{\bx}A)\ba|\, < \,\psi(\|\ba\|^n)\text{ for some }a_0\in \mathbb{Z}\}\,.\]
Observe that  $\,\,\limsup \mathcal {L}(\mathbf {a})$ is the projection of $\cW^{\theta}_{n}(\psi) \cap \mathscr{H} $ onto $\R^\xd$.
Here and elsewhere, unless stated otherwise, $\mathbf {a} \in \mathbb{Z}^n\setminus \{\mathbf{0}\}$ and any unspecified limsup is taken over such $\ba$. Since the projection from $\mathscr{H}$ to $\R^\xd$  is bi-Lipschitz,
Theorem~\ref{main convergence} will follow on showing that $|\limsup \mathcal {L}(\mathbf {a})| = 0$. In fact, it is sufficient to show that for each $\bx\in U$, we can choose  an open ball $B$ centred at $\bx$  with $11B\subseteq U$ such that
\begin{equation}\label{svhonest} |\limsup \mathcal {L}(\mathbf {a},B)|=0\,,
\end{equation} where
$\mathcal {L}(\mathbf {a},B):=\mathcal {L}(\mathbf {a})\cap B$.

To prove the above measure zero statement, it is natural  to consider
separately  the case when we have a `large derivative' and the case when we do not. More precisely,   we split the set  $\mathcal {L}(\mathbf {a},B)$ into two subsets depending on the size of the   quantity
$\nabla(\hat{\theta} (\bx)+ (\bx, \tilde{\bx}A)) =\nabla(\hat{\theta} (\bx))+[{\rm Id}_\xd\,\,A']\ba \, $ --  here  $A'$   is as introduced at the start of
\S\ref{sec:exp} and as usual $ \nabla $ denotes the gradient operator. With this in mind, for any sufficiently small open ball $B$ with $11B\subseteq U$, we define
\begin{equation}\mathcal {L}_{small}(\mathbf {a},B) = \left\{\bx \in \mathcal
{L}(\mathbf {a},B)~:~
\|\nabla\big(\hat{\theta} (\bx)+ (\bx, \tilde{\bx}A)\cdot \ba\big)\|<\sqrt{n\xd L
\|\ba\|}\right \}\end{equation}
where
\begin{equation}\label{eqn:L}
L:=\max \left \{\sup_{|\beta|=2 \, ,  \ \mathbf{x}\in
		2B}\|\partial_{\beta}\hat{\theta}(\bx)\|, \frac{1}{4r^2}\right \}
\end{equation}
and $r$ is the radius of $B$.
Here for a multi-index $\beta=(i_1,\dots,i_\xd)$ of non-negative integers $|\beta|:=i_1+\dots+i_\xd$ and $\partial_\beta$ denotes the corresponding differentiation operator, that is $\frac{\partial^{|\beta|}}{\partial x_1^{i_1}\dots\,\partial x_\xd^{i_\xd}}$. Set $\mathcal
{L}_{large}(\mathbf {a}, B) = \mathcal {L}(\mathbf {a},B)
\backslash \mathcal {L}_{small}(\mathbf {a},B)$.   We will prove that for any `appropriately chosen' $B\subseteq 11B\subseteq U$,

 \begin{equation}\label{Borcant1}
     |\limsup  \mathcal {L}_{large}(\mathbf {a},B)|=0 \end{equation}  \noindent and
 \begin{equation}\label{Borcant2}
 |\limsup  \mathcal {L}_{small}(\mathbf {a},B)|=0 \, . \end{equation}

\noindent Clearly, on combining the measure zero statements \eqref{Borcant1}  and \eqref{Borcant2}  we obtain the desired measure zero statement \eqref{svhonest}.

\section{Estimating the measure of $\displaystyle \limsup \mathcal
{L}_{large}(\mathbf {a},B)$} \label{poll}
In this section, we will establish (\ref{Borcant1}) as a simple consequence of the following statement.

\begin{proposition}\cite[Lemma 2.2]{BKM}\label{prop:large} Let $B\subseteq \mathbb{R}^\xd$
be a ball of radius $r$ and
$F \in C^2 (2B)$ where $2B$ is the ball with the same
centre as $B$ and radius $2r$.
Let \begin{equation}\label{defM*}
M^*:=\sup_{|\beta|=2,\,\mathbf{x}\in
2B}\|\partial_{\beta}F(\bx)\|
\end{equation}and
\begin{equation}\label{defM}
 M:=\max  \left \{M^*, \frac{1}{4r^2}\right \}.
\end{equation} Then for every $\delta'\, > \,0$, the set
of all $\mathbf{x}\in B$ such that $|p+
F(\bx)|\, < \,\delta'$ for some $p\in \mathbb{Z}$
and
\begin{equation}
\|\nabla F(\bx)\|\,\geq\,\sqrt{\xd M}
\end{equation}
\noindent has measure at most $K_{\xd }\delta'|B|$, where $K_{\xd }>0 $ is a constant dependent only on $\xd $.
\end{proposition}

To prove (\ref{Borcant1}) from Proposition \ref{prop:large}, we start with any sufficiently small open ball small $B$ in $U$. We fix $\ba \in \bbZ^n \setminus \{\mathbf{0}\}$ and
take $F(\bx)=((\mathbf{x},\tilde{\mathbf{x}}A), \hat{\theta}(\bx))\cdot (\ba,1)$
for $\bx \in 2B$ and $\delta'= \psi(\|\ba\|^n)$. Clearly
$M=L$, where $L$ is given by (\ref{eqn:L}). Hence by Proposition \ref{prop:large},
we get that
\[|\mathcal {L}_{large}(\mathbf {a},B)|\,\leq\,K_\xd \psi(\|\ba\|^n)|B|\,,\] and thus
\[\sum_{\ba \in \bbZ^{n}\setminus \{\mathbf{0}\}} |\mathcal {L}_{large}(\mathbf{a},B)| \ \leq \  K_\xd  \!
\sum_{\ba \in \bbZ^{n}\setminus \{\mathbf{0}\}} \psi(\|\ba\|^n) |B|\ll
\sum_{h=1}^\infty h^{n-1}\psi(h^n)\ll \sum_{h=1}^\infty \psi(h).\]
Since the latter sum is convergent, on applying the  Borel-Cantelli lemma we obtain (\ref{Borcant1}), as desired. In the above as well as elsewhere $\ll$ means an inequality with an unspecified multiplicative constant.

\vspace*{0ex} ~ $ \hfill \Box$      \\

In order to establish \eqref{Borcant2}   we will use the `inhomogeneous transference principle' introduced in
\cite[Section 5]{BeVe}, which simplified version is recalled in the next section.

\section{Inhomogeneous transference principle}\label{ITP}

Throughout this section, we shall let $B$ be an open ball in $\bbR^\xd$
and $\varrho$ be the $\xd $-dimensional Lebesgue measure restricted to $B$ so that the closed ball $\overline{B}$ becomes the support
of $\varrho$.

Consider two countable index sets $\mathcal{T}$, $\mathcal{A}$ and two maps $H: (t,\alpha,\eta)\mapsto H_{t}(\alpha,\eta)$ and
$I:(t,\alpha,\eta)\mapsto I_{t}(\alpha,\eta)$ from
 $\mathcal{T}\times \mathcal{A}\times \bbR_{+}$ to the collection of all open sets in $\bbR^\xd$. Take a set $\Phi$ of functions $\phi:\mathcal{T}
 \longrightarrow \bbR_{+}$. For each $\phi\in \Phi$, define
 \[\Lambda_I (\phi):= \displaystyle \limsup_{t\in \mathcal{T}} \bigcup_{\alpha\in \mathcal{A}}I_{t}(\alpha,\phi(t))  \quad   \text{ and }   \quad
 \Lambda_H (\phi):= \displaystyle \limsup_{t\in \mathcal{T}} \bigcup_{\alpha\in \mathcal{A}}H_{t}(\alpha,\phi(t))\,.\]We now discuss the two main
 properties which enables one to transfer zero $\varrho$-measure statements for the `homogeneous' limsup
 sets $\Lambda_H (\phi)$ to the `inhomogeneous' limsup sets $\Lambda_I (\phi)$.
 \begin{enumerate}
  \item \textbf{Intersection property}: The triplet $(H,I,\Phi)$ is said to satisfy the intersection property if for any
  $\phi \in \Phi$, there exists $\phi^{*}\in \Phi$ such that for all but finitely many $t\in \mathcal{T}$ and for all distinct $\alpha ,\alpha'\in \mathcal{A}$, we have
  \begin{equation}\label{eqn:intersection}
   I_t(\alpha, \phi(t))\cap I_t(\alpha',\phi(t))\subseteq \bigcup_{\alpha''\in \mathcal{A}}H_t(\alpha'',\phi^{*}(t))\,.
  \end{equation}
  \item \textbf{Contracting property}: We say that $\varrho$ is  contracting with respect to $(I,\phi)$ if for any $\phi\in \Phi$, there
  exists $\phi^{+}\in \Phi$, a sequence of positive numbers $\{k_t\}_{t\in \mathcal{T}}$ with $ \sum_{t\in
  \mathcal{T}}k_t\, <  \,\infty$ and for all but finitely many $t\in \mathcal{T}$ and all $\alpha \in \mathcal{A}$, a
  collection $\mathcal{C}_{t,\alpha}$
  of balls $\mathfrak{B}$ centred in $\overline{B}$ satisfying the three conditions given below:
  \begin{equation}\label{contract1}
  \overline{B}\cap I_{t}(\alpha, \phi(t))\subseteq \displaystyle \bigcup_{\mathfrak{B}\in \mathcal{C}_{t,\alpha}}\mathfrak{B}\,,
  \end{equation}
  \begin{equation}\label{contract2}
   \overline{B}\cap \displaystyle \bigcup_{\mathfrak{B}\in \mathcal{C}_{t,\alpha}}\mathfrak{B} \subseteq I_t(\alpha, \phi^{+}(t))
  \end{equation} and
  \begin{equation}\label{contract3}
   \varrho(5\mathfrak{B}\cap I_t(\alpha, \phi(t)))\leq k_t \varrho(5\mathfrak{B})\,.
  \end{equation}
  \end{enumerate}
  The main transference for our purpose, which follows easily from \cite[Theorem 5]{BeVe}, can be stated as follows.
  \begin{theorem}\label{thm:transference}
 If $(H,I,\Phi)$ satisfies the intersection property and $\varrho$ is contracting with respect to $(I,\phi)$ then
 \[\forall \phi \in \Phi\,\,\varrho(\Lambda_H(\phi))=0\qquad\Longrightarrow\qquad\forall \phi \in \Phi\,\,\varrho(\Lambda_I(\phi))=0\,.\]
  \end{theorem}

 \section{Proof of (\ref{Borcant2}) from Theorem \ref{thm:transference}} \label{svsv}
Let $\mathcal{T}:= \mathbb{Z}_{+}, \mathcal{A}=(\mathbb{Z}^n \setminus \{\mathbf{0}\})\times \mathbb{Z}$. For $t\in \mathcal{T},
\alpha:= (\mathbf{a},a_0)\in \mathcal{A}$ and $\eta \in \bbR_{+}$, we set

\begin{equation}\label{inhomset}
I_t(\alpha,\eta):= \displaystyle \left \{\bx\in U:
   \left|
\begin{array}{ll}
|\hat{\theta}(\bx)+ a_0 + (\bx, \tilde{\bx}A)\ba| < \displaystyle
\frac{\eta}{2^{nt}}\\[2ex]
\|\nabla(\hat{\theta} (\bx) +(\bx, \tilde{\bx}A) \cdot \ba) \| < \displaystyle
\sqrt{n\xd L}\times \eta \times 2^{t/2}\\[2ex]
    2^t \leq \|\ba\|< 2^{t+1}
\end{array}\right.\right\} \end{equation} and
\begin{equation}
H_t(\alpha,\eta):= \displaystyle \left \{\bx\in U:
   \left|
\begin{array}{ll}
| a_0 + (\bx, \tilde{\bx}A)\ba| < \displaystyle
\frac{2\eta}{2^{nt}}\\[2ex]
\|\nabla (\bx, \tilde{\bx}A) \cdot \ba \| < \displaystyle
2\sqrt{n\xd L}\times \eta \times 2^{t/2}\\[2ex]
    \|\ba\|< 2^{t+2}
\end{array}\right.\right\}\,. \end{equation}

\noindent Given $\delta \in \bbR$, define \[\phi_{\delta}: \mathcal{T}\longrightarrow \bbR_{+} : t \to  \phi_{\delta}(t):=2^{\delta t} \,.\]
Pick $\gamma \, >  \,0$ and consider the set
\[\Phi:= \{\phi_{\delta}: 0\leq \delta\, <  \,\gamma\}\,.\]

 Recall,  that in view of  \S\ref{poll},   the proof of Theorem~\ref{main convergence} has been reduced to showing  the truth of  (\ref{Borcant2}).   The above transference principle plays a key role in carrying out this task. With this in mind, the proof of (\ref{Borcant2}) splits naturally  into three main steps.
 Let $B\subseteq U$ be an open ball and  recall that  $\varrho$ is  the $\xd $-dimensional Lebesgue measure restricted to $B$.
\begin{enumerate}
\item[]
\begin{enumerate}
  \item[Step 1:]  We show that if $\gamma$ is appropriately chosen, then
 \begin{equation}\label{eqn:homo}
 \varrho(\Lambda_{H}(\phi_{\delta}))=0   \quad \forall    \ \delta \in [0,\gamma) \, ,
 \end{equation} regardless of the choice $B$.  This will be the subject of \S\ref{SVisGood} -- \S\ref{homo}. \\
 \item [Step 2:] We show that the triplet $(H,I,\Phi)$ as  defined above satisfies the intersection property.  This will be the subject of \S\ref{verinter}. \\
 \item[Step 3:]  We show that  $\varrho$ is contracting with respect to $(I,\phi_{\delta})$.  This will be the subject of \S\ref{vercont}.
\end{enumerate}
\end{enumerate}
The upshot  of successfully carry out these steps, is that on applying  Theorem~\ref{thm:transference}, we are able to conclude that
\[ \varrho(\Lambda_{I}(\phi_{\delta}))=0   \quad \forall    \ \delta \in [0,\gamma) \, .\]
This  in turn implies  (\ref{Borcant2}), since
\[\limsup \mathcal {L}_{small}(\mathbf {a},B)   \ \subseteq   \ \Lambda_{I}(\phi_{\delta}) \cap B \quad \forall   \ \delta   >  \,0\,. \]

 To carry out Step 1  we shall
employ dynamical tools. For that, we need to recall a
few elementary
properties of the so called `good functions' introduced by Kleinbock and  Margulis \cite{KM}.

\section{$(C, \alpha)$-good functions}  \label{SVisGood}

\noindent Let $C$ and $\al$ be positive numbers and $V$ be an open
subset of
$\bbR^\xd$. A function $f : V \to \bbR$ is said to be
$(C,\al)$-good on
$V$ if for any open ball $B \subseteq V$,~and
for any $\varepsilon > 0$, one has:

\begin{equation}\label{gooddef}
\bigg| \bigg\{ \bx \in B  :  |f(\bx)| < \varepsilon
 \bigg \} \bigg| \leq
C\left(\displaystyle \frac{\varepsilon}{\sup_{\bx \in
B}|f(\bx)|}\right)^{\al}|B|.
\end{equation} \\
Now consider  $\mathbf{f}=(f_1,\dots,f_n)$, a map from an open subset $U\subseteq \bbR^\xd$ to $\bbR^n$. We will say that $\mathbf{f}$ is good at $\bx_0\in U$ if there exists a neighbourhood $V\subseteq U$ and $C, \alpha>0$ such that any linear combination of $1,f_1,\dots,f_n$ is $(C,\alpha)-good$ on $V$. The map $\mathbf{f}$ is said to be good if it is good at every point of $U$. Note that $C,\alpha$ need not be uniform. \\

\noindent We will make use of the following properties of $(C,
\al)$-good functions, e.g. see \cite{Kleinbock-extremal}.
\begin{enumerate}
\item[(G1)] If $f$ is $(C,\al)$-good on an open set $V$, so
is $\lambda
f$ for all $\lambda \in
\bbR$.\\
\item[(G2)] If $f_i, i \in I$ are $(C,\al)$-good on $V$, so
is $\sup_{i \in
I}|f_i|$.\\
\item[(G3)] If $f$ is $(C,\al)$-good on $V$ and for some
$c_1,c_2\, >  \,0,\, c_1\leq
\frac{|f(\bx)|}{|g(\bx)|}\leq c_2
\text{ for all }\bx \in V$, then g is
$(C(c_2/c_1)^{\al},\al)$-good on $V$.\\
\item[(G4)] If $f$ is $(C,\al)$-good on $V$, it is
$(C',\alpha')$-good
on $V'$ for every $C' \geq \max\{C,1\}$, $\alpha' \leq \alpha$ and $V'\subset V$.
\end{enumerate}
One can note that from (G2), it follows that the supremum
norm of a vector valued function $\f$ is $(C,\al)$-good
whenever each of its components is $(C,\al)$-good.
Furthermore, in view of (G3), we can replace the norm by an
equivalent one, only affecting
$C$ but not $\al$. \\

The following result provides us with  an important class of
good functions.
\begin{proposition}\cite[Lemma 3.2]{BKM}\label{goodprop} Any polynomial $f\in
\mathbb{R}[x_1,...,x_\xd]$ of degree not
exceeding $l$ is $(C_{\xd ,l},\frac{1}{\xd l})$-good on
$\mathbb{R}^\xd$, where
$C_{\xd ,l}=\frac{2^{\xd +1}\xd l(l+1)^{1/l}}{V_\xd}$. In particular,
constant and linear polynomials are
$(\frac{2^{\xd +2}\xd}{V_\xd},\frac{1}{\xd })$-good on $\mathbb{R}^\xd$.\end{proposition}

The main dynamical tool that we will be exploiting to show (\ref{eqn:homo}) is commonly known as the ``quantitative nondivergence"
estimate in the space of unimodular lattices. This constitutes our next section.
\section{A quantitative nondivergence estimate}
Let $W$ be a finite dimensional real vector space. For a
discrete subgroup $\Gamma$ of $W$, we set $\Gamma_{\bbR}$ to
be the
minimal linear subspace of $W$ containing $\Gamma$. A
subgroup $\Ga$ of
$\La$ is said to be primitive in $\La$ if $\Ga = \Ga_{\bbR}
\cap
\La$. We denote the set of all nonzero primitive subgroups
of $\Gamma$ by $\mathcal{L}(\Gamma)$. Let $j:=\dim
(\Gamma_{\bbR})$ be the \emph{rank} of $\Gamma$. We say that
$\bw\in \bigwedge^j (W)$ represents $\Gamma$ if
\[\bw=\left \{\begin{array}{rcl}1 &\text{if
} & j=0\\[1ex] \bv_1\wedge\dots\wedge\bv_j &\text{if }&j>0\text{ and
}\bv_1, \dots,
\bv_j \text{ is a basis of }\Gamma\,.\end{array}\right.\]
In fact, one can easily see that such a representative of
$\Gamma$ is always unique up to a sign.

A function $\nu:\bigwedge (W)\longrightarrow \bbR_{+}$ is
called \emph{submultiplicative} if
\begin{itemize}
\item[]
\begin{itemize}
\item [(i)]$\nu$ is continuous with respect to natural
topology on $\bigwedge (W)$,
\vspace{2ex}
\item [(ii)]$\forall \ t\in \bbR \text{ and }\bw\in \bigwedge
(W)$, \ $\nu(t\bw)=|t|\nu(\bw)$,
\vspace{2ex}
\item [(iii)]$\forall \ \mathbf{u},\bw \in \bigwedge
(W),  \  \nu(\mathbf{u}\wedge\bw)\leq \nu(\mathbf{u})\nu(\bw)$.
\end{itemize}
\end{itemize}

\noindent In view of property (ii) above, without
any confusion,  we can define $\nu(\Gamma):=\nu(\bw) $ where  $\bw$ represents $\Gamma$.  Armed with the notion of submultiplicative,  we are in the position to  state the  ``quantitative nondivergence''
estimate that we will require in establishing (\ref{eqn:homo}).

\begin{theorem}\cite[Theorem 6.2]{BKM} \label{BKM1}
Let $W$ be a finite dimensional real vector space, $\Lambda$
a discrete subgroup of $W$ of rank $k$, and let a ball
$B = B(\bx_0,r_0) \subset
\bbR^\xd$  and a continuous map $H : \tilde{B} \rightarrow
\GL(W)$ be given, where $\tilde{B}:=B(\bx_0,3^kr_0)$. Take
$C\geq1$, $\alpha > 0,~0 < \rho < 1$ and let $\nu$ be a
submultiplicative function on $\bigwedge(W)$. Assume that for
any $\Gamma \in \mathcal{L}(\Lambda)$,
\begin{itemize} \item[]
\begin{itemize}
\item[(KM1)] the function $\bx \mapsto \nu(H(\bx)\Gamma)$ is
$(C,\alpha)$-good
on
$\tilde{B}$,
\vspace{2ex}
\item[(KM2)] $\displaystyle \sup_{\bx\in B}\nu(H(\bx)\Gamma)
\geq
\rho$,
\vspace{2ex}
\item[(KM3)] $\forall~\bx \in \tilde{B}$,   \ $\# \{ \Gamma \in
\mathcal{L}(\Lambda):\nu( H(\bx)\Gamma)
 < \rho\} < \infty$.
 \end{itemize}
\end{itemize}

\noindent Then for every  $\varepsilon'' >0$  we have that
\begin{equation}\label{BKMeq}
\Big|\big\{\bx \in B \, : \, \nu(H(\bx)\mathbf{\lambda  }) < \varepsilon'' ~
\text{for
some}~ \mathbf{\lambda} \in \La \backslash \{0\} \, \big\}\Big| \, < \, k \,(3^\xd
N_\xd)^{k} \,
C \left(\frac{\varepsilon''}{\rho}\right)^{\alpha}|B|.
\end{equation}
\end{theorem}

\section{Proof of (\ref{eqn:homo})}   \label{homo}

Fix a ball $B  \subset U$ such that  $11B  \subset U$. For $t\in \mathbb{Z}_+$ and $\delta \in [0,\gamma)$, we define
the set

\[\begin{array}{rl}\mathcal{A}_t :=  \displaystyle
\bigcup_{\alpha\in \mathcal{A}} \ (\,H_{t}\big(\alpha,\phi_{\delta}(t)) \, \cap  \,  B  \, \big)
&= \left \{\rule{0ex}{8ex}\bx\in B:\exists (\mathbf{a},a_0)\in \bbZ^{n}\setminus
\{\mathbf{0}\}\times \bbZ \text{
s.t.}\right.\\[1ex]
 & \hspace{5ex} \left|\rule{0ex}{6ex}\right.
\left.\begin{array}{ll}
|a_0 + (\bx, \tilde{\bx}A)\ba| < \displaystyle
\frac{2\times 2^{\delta t}}{2^{nt}}\\[1ex]
\|\nabla (\bx, \tilde{\bx}A) \cdot \ba \| < \displaystyle
2\sqrt{n\xd L}\times 2^{\delta t}\times 2^{t/2}\\[1ex]
    \|\ba\|< 2^{t+2}
\end{array}\rule{0ex}{8ex}\right\}\,.\end{array}\]

\noindent Then, by definition
 $$
  \ \Lambda_{H}(\phi_{\delta}) \cap B   \ \subseteq \  \textstyle{\limsup_{t \to \infty}} \mathcal{A}_t
 $$
 and so \eqref{eqn:homo} follows on showing that
\begin{equation}\label{svmid}
 |  \textstyle{\limsup_{t \to \infty}}  \mathcal{A}_t  | = 0 \ .
\end{equation}

\noindent With this in mind, pick $\beta\in \big(0,\frac{1}{2(n+1)}\big)$ and set
\begin{equation}\label{eqn:constants1}
\delta':=\displaystyle \frac{2}{2^{nt}},  \qquad
K:=\displaystyle2\times \sqrt{n\xd L}\times 2^{t/2}, \qquad  T:=2^{t+2},
\end{equation}
\begin{equation}\label{eqn:constants2a}
\varepsilon':=(\delta' K
T^{n-1})^{\frac{1}{n+1}} \ = \
\left(2^{2n}\sqrt{n\xd L}\right)^{\frac{1}{n+1}}
\frac{1}{2^{t/2(n+1)}}
 \end{equation}
 and
 \begin{equation}\label{eqn:constants2}
\varepsilon:= 2^{\beta t}\varepsilon'  \ =  \
\left(2^{2n}\sqrt{n\xd L}\right)^{\frac{1}{n+1}}
\frac{2^{\beta t}}{2^{t/2(n+1)}}.
  \end{equation}

\noindent Furthermore, for $ \bx \in \R^\xd$,  let
\begin{equation}\label{unip} u_{\mathbf{x}} : = \left
(\begin{array}{rccl}
1 & 0 & \bx& \bx A'+\mathbf{a}_0\\
0 & I_\xd & I_\xd& A'\\
0 & 0 & I_{n}\\
\end{array}\right)\end{equation}
\noindent and for  $t\in \mathbb{Z}_+$, let
\begin{equation}\label{diag}
g_t := \diag
\left(\frac{\varepsilon}{\delta'},\frac{\varepsilon}{K},\dots,\frac{\varepsilon}{K},
\frac{\varepsilon}{T},\dots,\frac{\varepsilon}{T}\right)\,,
\end{equation}
where $\varepsilon,\delta',T,K$ are defined above. Note that these parameters depend on $t$ and some fixed constants.
Also, denote by $\La$ the subgroup of $\bbZ^{1+\xd +n}$
consisting of vectors of the form:
\begin{equation}\label{latt}
\Lambda = \left\{\begin{pmatrix}
p\\
0\\
\vdots\\
0\\
\bq
\end{pmatrix} :  p \in \bbZ, \bq \in \bbZ^{n} \right\}.
\end{equation}

\noindent Then, it readily follows from the above definitions that

\begin{equation}\label{subset}
\mathcal{A}_t    \ \subseteq  \   \tilde{\mathcal{A}}_t:=\left\{\bx \in B
~:~ \|g_tu_\bx \lambda\| < 2^{\delta t} \varepsilon
~\text{for some}~\lambda \in \La \backslash \{0\}\right\} ,
\end{equation}

\noindent and so \eqref{svmid} follows on showing that

\begin{equation*}\label{svvmid}
 |  \textstyle{\limsup_{t \to \infty}}  \tilde{\mathcal{A}}_t  | = 0 \ .
\end{equation*}

\noindent In view of the Borel-Cantelli lemma,   this will follow on showing that

\begin{equation}\label{eqn:Borecant3}
\displaystyle \sum_{t=0}^{\infty} |\tilde{\mathcal{A}}_t|<\infty\,.
\end{equation}

With the intention of using Theorem \ref{BKM1} to prove
(\ref{eqn:Borecant3}), we take
$W=\bbR^{1+\xd +n}$ with basis
$\mathbf{e}_0,\mathbf{e}_{*1},\dots,\mathbf{e}_{*\xd},\mathbf{e}_1,\dots,\mathbf{e}_n$,
$\Lambda$ as
given by  (\ref{latt}) and $H(\bx)=g_tu_{\bx}$. The
submultiplicative function $\nu$ on $W$ is  chosen as
described in \cite[\S7]{BKM}.  Namely, let $W_*$ be the $\xd $-dimensional subspace of $W$ spanned by $\esa, \dots,
\ess$ so that $\Lambda$ given by \eqref{latt} is equal to the intersection of $\bbZ^{1+\xd +n}$ and $W_{*}^{\perp}$. Here we  identify $W_{*}^{\perp}$ with $\bbR^{n+1}$
canonically. Also let $\cW$ be the
ideal of $\bigwedge(W)$ generated by $\bigwedge^{2}(W_*)$ and let
$\pi_{*}$ be the orthogonal
projection with kernel $\cW$.  Then $\|\bw\|_e$ is defined to  be the
Euclidean norm of $\pi_{*}(\bw)$. In other  words,
if $\bw$ is written as a sum of exterior products of the
base vectors $\mathbf{e}_i$ and $\mathbf{e}_{*i}$,
to compute $\nu(\bw)$ we ignore the components containing
exterior products of the type
$\mathbf{e}_{*i}\wedge\mathbf{e}_{*j},1\leq i\neq j\leq \xd$,
and
simply take the Euclidean norm of the sum of the remaining components. By definition, it is immediate that
$\nu|_W$ agrees with the Euclidean norm. \\

For appropriately  determined quantities  $C,\al,\rho$  we now validate, one by one, the conditions (KM1)-(KM3) associated with Theorem \ref{BKM1}. Condition (KM3) can be verified  for any
$\rho\leq1$ in exactly the same manner as in
\cite[\S7]{BKM}. For the
verification of the remaining conditions, we begin with the explicit computation of the quantity
$H(\bx)\bw$ for any $\bw\in \bigwedge^k(W_{*}^{\perp})$ and
$k=1,\dots,n+1$. On  writing
$\bx=(x_1,\dots,x_\xd)$ and
$(\bx,\tilde{\bx}A)=(f_1(\bx),\dots,f_n(\bx))$, we see that
\begin{enumerate}
\item $H(\bx)\,\mathbf{e}_0= \frac
{\varepsilon}{\delta'}\,\mathbf{e}_0$,
\vspace{1ex}
\item $H(\bx)\,\mathbf{e}_{*i}= \frac
{\varepsilon}{K}\,\mathbf{e}_{*i}\
 \text{  \ for  all \ }1\leq i\leq \xd$,
\vspace{2ex}
\item $H(\bx)\,\mathbf{e}_{i}= \frac
{\varepsilon}{\delta'}\,f_i(\bx)\,\eo + \frac
{\varepsilon}{K}\,\sum_{j =
1}^{\xd } \frac{\partial f_j (\bx)}{\partial x_i} \esj + \frac {\varepsilon}{T}\,\ei\, \text{ \ for all \
}1\leq i\leq n.$
 \end{enumerate}
Note that each $f_i(\bx)$ is a polynomial in $x_1,\dots,x_\xd$
with degree at most $1$ so that each  partial derivative $\frac{\partial f_j (\bx)}{\partial x_i}$ is constant.

\subsection{Checking (KM1)}
Since $\Lambda=\bbZ^{1+\xd +n}\cap
W_{*}^{\perp}$, any representative $\bw\in \bigwedge^k(W)$
of any subgroup of $\Lambda$ of rank $k$, $1\leq k\leq n+1$,
can be written as $
\sum_{I}a_{I}\mathbf{e}_{I}$, where each $a_I\in \bbZ$ and
$\mathbf{e}_I=\mathbf{e}_{i_i}\wedge \dots \wedge
\mathbf{e}_{i_k}$ with $i_1,\dots,i_k \in \{0,1,\dots,n\},
i_1<\cdots<i_k$. \\

Since each component of $\pi_{*}(H(\bx)\bw)$ is a polynomial
in $x_1,\dots,x_\xd$
 with degree at most $1$ and in view of Proposition \ref{goodprop}, each of them is
$(\frac{2^{\xd +2}\xd}{V_\xd},\frac{1}{\xd })$-good on $\tilde{B}$.
This implies that the function
$\bx\mapsto\|\pi_{*}(H(\bx)\bw)\|$  is
$(\frac{2^{\xd +2}\xd}{V_\xd},\frac{1}{\xd })$-good on $\tilde{B}$. As

\[\displaystyle \frac{1}{2^{\frac{1+\xd +n}{2}}} \ \leq  \
\frac{\|\pi_{*}(H(\bx)\bw)\|}{\nu(\pi_{*}(H(\bx)\bw))} \ \leq \
1\,,\]
it follows from property
(G4) of good functions that  $\nu(\pi_{*}(H(\bx)\bw))$ is
$(C,\al)$-good on $\tilde{B}$ with
 \begin{equation}\label{eqn:C,al}
C:=\max
\left\{\frac{2^{\left(\xd+2+\frac{1+\xd +n}{2\xd}\right)}\xd}{V_\xd},1\right\}\text{
 \ and  \  }\al:=\frac{1}{\xd}\,.
 \end{equation} This verifies condition (KM1).

 \subsection{Checking (KM2)}\label{ckm2}
Let $\Gamma$ be a subgroup of $\Lambda$ with rank $k$ and
$\bw\in \bigwedge^k( W_{*}^{\perp})$ represent $\Gamma$. We
first consider the case $k=n+1$.
Thus,  $\bw = w\, \eo \wedge \ea \wedge \dots \wedge \en$,  where
$w\in \bbZ\backslash\{0\}$. Hence, for any $\bx\in B$,  it is easily verified  that the
coefficient of
$\eo \wedge \esa \wedge \mathbf{e}_2 \wedge \dots \wedge
\en$ in $\pi_{*}(H(\bx)\bw)$ is
\[w\frac{\varepsilon^{n+1}}{\delta' K T^{n-1}}\,.\] \\
It now follows via  \eqref{eqn:constants1}--(\ref{eqn:constants2}),  that
\begin{eqnarray}\label{eqn:toprank}
\displaystyle \sup_{\bx\in
B}\nu\big(H(\bx)\Gamma\big)=\displaystyle \sup_{\bx\in
B}\nu\big(H(\bx)\bw\big)  & \geq & \sup_{\bx\in B}
\|\pi_{*}(H(\bx)\bw)\|  \nonumber \\[2ex]   & \geq   &
\Big|w\frac{\varepsilon^{n+1}}{\delta' K T^{n-1}}\Big|\nonumber
\\[2ex]  &  =   & |w|2^{\beta(n+1)
t}\frac{(\varepsilon')^{n+1}}{\delta' K T^{n-1}}\nonumber
\\[2ex]
& = & |w|2^{\beta(n+1) t} \geq
1.
\end{eqnarray}
Thus, when $k=n+1$ condition (KM2) is valid for any $0< \rho  < 1$.  \\

Assume now that $1\leq k\leq n$. To bound the norm of
$\|\pi_{*}(H(\bx)\bw)\|$ from below, we will proceed along
the lines of \cite[\S5.3]{G-Monat}
using a technique from \cite{Kleinbock-exponent}. As
observed in \cite[\S5.3]{G-Monat}, for any $\bx\in
B $
$$\|\pi_{*}(H(\bx)\bw)\|\geq\|\tilde{g}_{t}\tilde{u}_{\bx}\bw\|$$
where
\begin{equation}\label{unipnew} \tilde{u}_{\bx} =
\begin{pmatrix}
1  & \bx&\tilde{\bx}A\\
0  & I_{n}\\
\end{pmatrix}   \quad  \text{ and  }  \quad
\tilde{g}_{t} =
\diag \left(\frac{\varepsilon}{\delta'},
\frac{\varepsilon}{T},\dots,\frac{\varepsilon}{T}\right). \\
\end{equation}
Hence
\begin{eqnarray}\label{eqn:toprank*}
\displaystyle \sup_{\bx\in
B}\nu\big(H(\bx)\Gamma\big)=\displaystyle \sup_{\bx\in
B}\nu\big(H(\bx)\bw\big) & \geq & \sup_{\bx\in B}
\|\pi_{*}(H(\bx)\bw)\|  \nonumber \\[2ex]   & \geq   &  \sup_{\bx\in B}
\|\tilde{g}_{t}\tilde{u}_{\bx}\bw\|   \, .
\end{eqnarray}

\noindent Thus, the name of the game is to  bound $\sup_{\bx\in
B}\|\tilde{g}_{t} \tilde{u}_{\bx} \bw\|$ from below.  It follows from ($4.6$) in \cite{Kleinbock-exponent}, that

\begin{equation}\label{boundformat1}
\displaystyle \sup_{\bx\in B}\|\tilde{g}_{t} \tilde{u}_{\bx}
\bw\|  \ \geq  \  \frac{1}{2^\frac{n+1}{2}}\max
\left\{\left(\frac{\varepsilon^k}{\delta'
T^{k-1}}\right)\sup_{\bx \in B}\|(\bx,\tilde{\bx}A) \
\bc(\bw)\|, \
\left(\frac{\varepsilon}{T}\right)^k\|\pi(\bw)\|\right\}
\end{equation} \\

\noindent where $ \bc $ is the function given by \eqref{defc},  $\pi$ is the projection from
$\bigwedge(W_{*}^{\perp})$ to $\bigwedge(W_{1\to n})$ and
$W_{1\to n}$ stands for the
span of the vectors  $\be_1,\dots,\be_n$ . Now, recall
  that
\begin{equation*}
(\bx,\tilde{\bx}A) = \tilde{\bx}R_A
\end{equation*}
\noindent where $R_A$ is given by (\ref{defR}). Therefore, we can replace $\sup_{\bx \in B}\|(\bx,\tilde{\bx}A)
\bc(\bw)\|$ in the above  norm
calculation  by $\sup_{\bx \in B} \|\tilde{\bx}R_A
\bc(\bw)\|$. As the functions $1,x_1,\dots,x_\xd$ are
linearly independent over $\bbR$ on
$B$, the map
$$\bv \mapsto \displaystyle \sup_{\bx\in
B}\|\tilde{\bx}\bv\|$$ defines a norm on $(\bigwedge (W_{1\to
n}))^{\xd +1}$ which must be equivalent to the supremum
norm on $(\bigwedge (W_{1\to n}))^{\xd +1}$. Thus, there is a
constant $K_2>0$ depending on $\xd ,n$ and $B$, such that
\[\sup_{\bx \in B} \|\tilde{\bx}R_A  \ \bc(\bw)\| \ \geq  \ K_2
\|R_{A} \ \bc(\bw)\|\,,\] \\
and consequently, via \eqref{boundformat1},  that

\begin{equation}\label{boundformat2}
\displaystyle \sup_{\bx\in B}\|\tilde{g}_{t} \tilde{u}_{\bx}
\bw\|  \ \geq  \  \frac{1}{2^\frac{n+1}{2}}\max
\left\{\left(\frac{\varepsilon^k}{\delta' T^{k-1}}\right)K_2
\|R_{A} \ \bc(\bw)\|,
\left(\frac{\varepsilon}{T}\right)^k\|\pi(\bw)\|\right\}\,.
\end{equation}\\

To continue,  we consider two separate cases depending on the size of the rank $k$. We first note that from Lemma $5.1$ in
\cite{Kleinbock-exponent},
we get that for any $n-\xd  < k \leq n$ for all but
finitely many $\bw \in \bigwedge^{k}(\Lambda)$ we have that
$\|R_A \ \bc(\bw)\| \geq 1$. Also note that $\|R_A \ \bc(\bw)\|$ does not vanish, as otherwise if $\|R_A \ \bc(\bw_0)\|$ were zero, then, by the linearity of the map $\bc(\bw)$, we would get $\|R_A \ \bc(\lambda\bw)\|=0$ for all integers $\lambda$, contrary to what we have already seen. Consequently, there is  a constant $K_3>0$ depending only  on $A$, such that
\begin{equation}
\|R_A \ \bc(\bw)\| \geq K_3\,.
\end{equation}
Therefore, by \eqref{boundformat2}, we get that
\begin{equation} \label{slv2}
\displaystyle \sup_{\bx\in B}\|\tilde{g}_{t} \tilde{u}_{\bx}
\bw\|\geq \frac{K_2 K_3}{2^\frac{n+1}{2}}\left
(\frac{\varepsilon^k}{\delta' T^{k-1}}\right).
\end{equation}\\
It follows from  \eqref{eqn:constants1}--(\ref{eqn:constants2}), that
\begin{eqnarray*}\displaystyle
\frac{\varepsilon^k}{\delta' T^{k-1}}  & =  &
\left(2^{2n}\sqrt{n\xd L}\right)^{\frac{k}{n+1}}\
\frac{1}{2^{\left(\frac{1}{2(n+1)}-\beta \right)kt}}\
\frac{2^{nt}}{2}\ \frac{1}{2^{(t+2)(k-1)}}\\[2ex]   & \geq   &
\displaystyle
\min_{n-\xd <k\leq n}\left(2^{2n}\sqrt{n\xd L}\right)^{\frac{k}{n+1}} \
\frac{1}{2^{\left(\frac{1}{2(n+1)}-\beta \right)nt}}\
\frac{2^{nt}}{2}\  \frac{1}{2^{(t+2)(n-1)}}\\[2ex]  & =  &
\displaystyle
\frac{1}{2^{2n-1}}  \  \ \min_{n-\xd <k\leq n}\left(2^{2n}\sqrt{n\xd L}\right)^{\frac{k}{n+1}}\
2^{\left(1-\left(\frac{1}{2(n+1)}-\beta\right)n\right)t}\,.
\end{eqnarray*}
On picking $\beta$ such that
\begin{equation}\label{svbeta1}
\frac{1}{2(n+1)} - \frac{1}{n}  \ <  \ \beta     \ < \  \frac{1}{2(n+1)} \, ,
\end{equation} \\
we obtain via \eqref{slv2},  that  for all subgroups $\Gamma$ of
$\Lambda$ with rank $n-\xd  < k \leq n$

\begin{equation}\label{eqn:lb1}
\displaystyle \sup_{\bx\in B}\|\tilde{g}_{t} \tilde{u}_{\bx}
\bw\|  \ \geq  \  \frac{K_2 K_3}{2^{\frac{5n-1}{2}}}\  \  \min_{n-\xd <k\leq n}\left(2^{2n}\sqrt{n\xd L}\right)^{\frac{k}{n+1}}  \, .
\end{equation} \\

\noindent We now obtain an  analogous lower bound result
for subgroups $\Gamma$ of
$\Lambda$ with rank  $1\leq k\leq n-\xd $.  In this case, a  consequence of (\ref{diocond})
is that there exist constants $0<\theta', K_4<1$, depending only on
$A$, such for any  $\bw\in
\bigwedge^k(\Lambda)$

\begin{equation}\label{diocond1}
\|R_{A} \
\bc(\bw)\|
 \ \geq  \ K_4 \,\|\pb(\bw)\|^{-\frac {(n - \theta')+1-{k}}{k}}.
\end{equation} \\
\noindent Also,   $\|\pi(\bw)\|\geq\|\pb(\bw)\|$ if $1\leq k\leq n-\xd $. Therefore, it follows from (\ref{boundformat2}) that

\begin{eqnarray}\label{boundformat3}
\displaystyle \sup_{\bx\in B}\|\tilde{g}_{t} \tilde{u}_{\bx}
\bw\|  &  \geq  &  \max
\left\{\left(\frac{\varepsilon^k}{\delta' T^{k-1}}\right)K_2
K_4 \,\|\pb(\bw)\|^{-\frac {(n - \theta')+1-{k}}{k}} \, ,  \
\left(\frac{\varepsilon}{T}\right)^k\|\pb(\bw)\|\right\}   \nonumber \\[2ex] & \geq & \  \kappa  \ \left(\frac{\varepsilon}{T}\right)^k  \, ,
\end{eqnarray}
where $\kappa$ is the solution to the equation
\begin{equation}
\frac{K_2K_4 T }{\delta'
}\,y^{-\frac{(n-\theta')+1-k}{k}}  \ =\  \,y  \, .
\end{equation}
\noindent In other words,  $\kappa := (K_2K_4)^{\frac{k}{n-\theta'
+1}}\left(\frac{T}{\delta'}\right)^{\frac{k}{n-\theta'+1}}$
and so it follows from \eqref{eqn:constants1}--(\ref{eqn:constants2}) that
\begin{eqnarray*}
 \kappa   \left(\frac{\varepsilon}{T}\right)^k \!\!\!  & = & \!\!\! (K_2K_4)^{\frac{k}{n-\theta'
+1}}\left(\frac{T}{\delta'}\right)^{\frac{k}{n-\theta'+1}}
\left(\frac{\varepsilon}{T}\right)^k  \\[2ex]
& =  & \!\!\! (K_2K_4)^{\frac{k}{n-\theta'
+1}}2^{\frac{k}{n-\theta'+1}}2^{\frac{(n+1)kt}{n-\theta'+1}}
\left(2^{2n}\sqrt{n\xd L}\right)^{\frac{k}{n+1}}
\frac{1}{2^{\left(\frac{1}{2(n+1)}-\beta \right)kt}}
\frac{1}{2^{(t+2)k}}  \\[2ex]
& =  &  \!\!\!
(K_2K_4)^{\frac{k}{n-\theta'
+1}}2^{\frac{k}{n-\theta'+1}}\frac{1}{2^{2k}}
\left(2^{2n}\sqrt{n\xd L}\right)^{\frac{k}{n+1}}
2^{\left(\left(\frac{n+1}{n-\theta'+1}-1\right)-
\left(\frac{1}{2(n+1)}-\beta \right)\right) kt}  \\[2ex]
& = &  \!\!\!
(K_2K_4)^{\frac{k}{n-\theta'
+1}}\left(2^{2n}\sqrt{n\xd L}\right)^{\frac{k}{n+1}}
2^{\left(\frac{1}{n-\theta'+1}-2\right)k}
2^{\left(\left(\frac{n+1}{n-\theta'+1}-1\right)-
\left(\frac{1}{2(n+1)}-\beta \right)\right) kt}.\end{eqnarray*}
On redefining  $ \beta $ if necessary; namely so that both  \eqref{svbeta1}  and \begin{equation}\label{svbeta2}
\frac{1}{2(n+1)} - \frac{\theta'}{(n- \theta')+1}  \ <  \ \beta     \ < \  \frac{1}{2(n+1)} \, ,
\end{equation} \\
hold,  it follows that
\begin{eqnarray*}
 \kappa   \left(\frac{\varepsilon}{T}\right)^k   & \geq  & K_5:= \min_{1\leq k\leq n-\xd }\displaystyle
(K_2K_4)^{\frac{k}{n-\theta'
+1}}\left(2^{2n}\sqrt{n\xd L}\right)^{\frac{k}{n+1}}
2^{\left(\frac{1}{n-\theta'+1}-2\right)k}\, .
 \end{eqnarray*}
Note that \eqref{svbeta2} has a solution $\beta$ since $0<\theta'<1$.
This together with \eqref{boundformat3} implies that   for all subgroups $\Gamma$ of
$\Lambda$ with rank  $1 \leq k \leq n-\xd $,   we have that
\begin{equation}\label{lb3}
\displaystyle \sup_{\bx\in B}\|\tilde{g}_{t} \tilde{u}_{\bx}
\bw\|  \ \geq   \ \frac{1}{2^\frac{n+1}{2}}\,K_5\,.
\end{equation}
On combining  (\ref{eqn:toprank}), \eqref{eqn:toprank*},
(\ref{eqn:lb1}) and (\ref{lb3}), we have verified condition (KM2)
with
\begin{equation}\label{eqn:rho}
\rho \ := \  \displaystyle \min \left \{\frac{1}{2},  \ \ \frac{K_2 K_3}{2^{\frac{5n-1}{2}}}\ \   \min_{n-\xd <k\leq n}
\left(2^{2n}\sqrt{n\xd L}\right)^{\frac{k}{n+1}} \!, \ \
 \frac{1}{2^\frac{n+1}{2}}\,K_5\right\}\,.
\end{equation} \\

We are now in the position  to apply Theorem \ref{BKM1} to establish the desired convergent sum statement
(\ref{eqn:Borecant3}).
\subsection{The proof of (\ref{eqn:Borecant3})} \label{homo3}
 With the choice of $\beta \in (0,1/2)$ made in the \S\ref{ckm2}, let
\begin{equation}\label{eqn:gamma}
0<\gamma<\frac{1}{2(n+1)}-\beta\, .
\end{equation}
Clearly,  $ \gamma>0$ and note that for any  $\delta\in [0,\gamma)$

\[\tilde{\mathcal{A}}_t\subseteq \Big\{  \bx \in B ~:~
\nu(H(\bx)\mathbf{\lambda}) < \sqrt{1+\xd +n}\,2^{\delta t}\,\varepsilon ~
\text{for
some}~ \mathbf{\lambda} \in \La \backslash \{0\} \Big\} \,.\] \\
Here we make use of the fact that  $\nu|_W$ coincides with the Euclidean norm on
$W$. Now on applying Theorem
\ref{BKM1} with $\varepsilon'':=\sqrt{1+\xd +n}\,2^{\delta t} \varepsilon $ where $ \varepsilon$ is given by  (\ref{eqn:constants2}),  and  $C,\al \text{ and }\rho$ are  as given
in (\ref{eqn:C,al}) and (\ref{eqn:rho}), we have

\begin{eqnarray} \label{ohya}
|\tilde{\mathcal{A}}_t| & \le   &  \Big|\big\{\bx \in B ~:~
\nu(H(\bx)\mathbf{\lambda}) < \sqrt{1+\xd +n}\,2^{\delta t}\,\varepsilon  \  \  ~
\text{for
some}~ \mathbf{\lambda} \in \La \backslash \{0\}  \; \big\}\Big| \nonumber \\[2ex]
& \leq  &  (n+1)(3^\xd N_\xd)^{n+1}
C
(1+\xd +n)^{\frac{1}{2\xd}}\left(\frac{\varepsilon}{\rho}\right)^{\frac{1}{\xd}}|B|  \nonumber \\[2ex]
& \leq   &  (n+1)(3^\xd N_\xd)^{n+1}
C
(1+\xd +n)^{\frac{1}{2\xd}}\frac{1}{\rho^{\frac{1}{\xd}}}\left(2^{2n}\sqrt{n\xd L}\right)^{\frac{1}{\xd(n+1)}}
\frac{1}{2^{\left(\frac{\frac{1}{2(n+1)}-(\beta+\delta)}{\xd}\right)t}}|B| \, .  \nonumber \\
& &\label{vb4}
\end{eqnarray}
As $\delta<\gamma$, it follows via  \eqref{eqn:gamma} that $ \delta+\beta < \frac{1}{2(n+1)}$,  and so
\[\displaystyle  \sum_{t=0}^{\infty}   |\tilde{\mathcal{A}}_t|  \ \ll \
\displaystyle  \sum_{t=0}^{\infty} 2^{- \left(\frac{\frac{1}{2(n+1)}-(\beta+\delta)}{\xd }\right)t} \ < \ \infty\,.\]
By the Borel-Cantelli Lemma, this establishes  (\ref{eqn:Borecant3}), as desired.

\section{Verification of the intersection property for $(H,I,\Phi)$} \label{verinter}
Let  $\gamma$  be given by  (\ref{eqn:gamma})  and let $\delta\in [0,\gamma)$. Suppose $t\in \mathcal{T}:= \mathbb{Z}_{+}$ is such
that $t(n-\delta)\geq 1$ and $\alpha:=(\ba, a_0), \alpha':=(\ba',a_0') \in \mathcal{A} $ with $\alpha\neq  \alpha'$.  Recall, $\mathcal{A} := (\mathbb{Z}^n \setminus \{\mathbf{0}\})\times \mathbb{Z} $.
Then, for any $\bx \in I_t(\alpha, \phi_{\delta}(t))\cap I_t(\alpha',\phi_{\delta}(t))$,  it is easily verified that

\begin{equation}\label{eqn:checkint}
	\left \{\begin{array}{rcl}|(a_0-a_0')+(\bx, \tilde{\bx}A)(\ba-\ba')| \ <  \ \displaystyle
	\frac{2\times 2^{\delta t}}{2^{nt}}\\[3ex] \|\nabla (\bx, \tilde{\bx}A) \cdot (\ba-\ba' )\| \ < \ \displaystyle
	2 \sqrt{n\xd L}\times 2^{\delta t}\times  2^{t/2}\\[2ex]
	\|(\ba-\ba')\| \ <  \  2^{t+2}\end{array}\right.\,.
\end{equation} \\
Suppose for the moment that  $\ba=\ba'$.   Then $a_0\neq a_0'$ as $\alpha\neq \alpha'$. This implies, in view of the first inequality  of (\ref{eqn:checkint}), that
$1\leq |(a_0-a_0')|<\frac{1}{2^{t(n-\delta)-1}}\leq 1$, which is a contraction. Thus $\ba\neq \ba'$  and so  $(\ba-\ba', a_0-a_0')\in \mathcal{A}$. The upshot of this together with $(\ref{eqn:checkint})$ is that $\bx\in H_t(\alpha'', \phi_{\delta} (t))$ with $\alpha''= (\ba-\ba', a_0 -a_0')$. This establishes $(\ref{eqn:intersection})$ with $\phi^{*}=\phi=\phi_{\delta}$ and thereby verifies the desired  intersection property associated with the Inhomogeneous Transference Principle.

\section{Verification of the contraction property of $\varrho$}  \label{section: contraction}  \label{vercont}

With reference to \S\ref{svsv}, recall that showing  $\varrho$ is contracting with respect to $(I,\phi_{\delta})$ is the third and final step in establishing Theorem~\ref{main convergence}.   We start by observing that in view of  \cite[Corollary 3.3]{Kleinbock-extremal} and the fact that the inhomogeneous function $\theta$ restricted to $\mathscr{H}$ is analytic\footnote{It is worth pointing out that this is the only  point in the proof of Theorem~\ref{main convergence}   where we  use the fact that $\theta|_{\mathscr{H}}$ is analytic.}, the functions $$\bx \mapsto |\hat{\theta}(\bx)+ a_0 + (\bx, \tilde{\bx}A)\ba|$$ and $$\bx\mapsto \|\nabla (\hat{\theta}(\bx)+(\bx, \tilde{\bx}A) \cdot \ba )\|$$ \\
defined on $U$ are good at every point of $U$. Now pick a point  $\bx_0 \in U$. On using property (G4) of good functions if necessary, we can choose an open ball $B$ with centre at $\bx_0$ and two positive constants $\mathfrak{C}, \mathfrak{\alpha_0}$ such that the above two functions are $(\mathfrak{C}, \mathfrak{\alpha_0})-good$ on $11B\subseteq U$. Throughout this section we fix such a ball $B$.

For each $t \in \mathcal{T} \text{ and }\alpha \in \mathcal{A}$,
consider the function $\mathbf{F}_{t,\alpha}:U \rightarrow \bbR$ given by
%

\[\mathbf{F}_{t,\alpha}(\bx)  := \max\left \{\begin{array}{rcl}2^{nt}  \ \sqrt{n\xd L} \ 2^{t/2}   \ |\hat{\theta}(\bx)+ a_0 + (\bx, \tilde{\bx}A)\ba| \, ,~\hspace*{8ex} \\[2ex]
 \|\nabla (\hat{\theta}(\bx)+(\bx, \tilde{\bx}A) \cdot \ba )\| \end{array}\right\}\,,\]

 \noindent where $L$ is given by (\ref{eqn:L}). It follows at once, from  the properties  of good functions, that for each $t \in \mathcal{T} \text{ and }\alpha \in \mathcal{A}$ we have that
\begin{equation}\label{goodeq} \mathbf{F}_{t,\alpha} \text{ is } (\mathfrak{C}, \alpha_0)-good \text{ on }  11B\,.
\end{equation}

\noindent Next, observe that for any  $\eta \in \bbR_{+}$ the first two inequalities appearing in  (\ref{inhomset})
are equivalent to the following single
inequality
\[\mathbf{F}_{t,\alpha}(\bx)<\eta \ \sqrt{n\xd L}\ 2^{t/2}\,.\] Hence,  for any  $t \in \mathcal{T}, \alpha=(\ba,a_0)\in \mathcal{A}
\text{ and }
\eta \in \bbR_{+}$
\begin{equation}\label{inhomsetmodi}
 I_{t}(\alpha, \eta)=\Big\{\bx\in U: \mathbf{F}_{t,\alpha}(\bx)<\eta \ \sqrt{n\xd L}\ 2^{t/2}\Big\}   \quad   \text{ if } \quad 2^t\leq \|\ba\|<2^{t+1},
\end{equation}

\noindent and
$ I_{t}(\alpha, \eta) = \emptyset$ otherwise.  For any $\delta \in [0,\gamma)$, consider the function  $\phi_{\delta}^{+} :\mathcal{T} \rightarrow \bbR_{+}$   given by
 \[\phi_{\delta}^{+}(t):= 2^{\frac{\delta+\gamma}{2} t}\  \,.\]\noindent Clearly,  $\phi_{\delta}^{+}\in \Phi$.   Also, for
any $t\in \mathcal{T}$ we have that
\begin{equation} \label{slv1} I_t (\alpha, \phi_{\delta} (t))  \; \subseteq  \; I_t (\alpha, \phi_{\delta}^{+} (t))  \, .
\end{equation}

In order to establish the desired contracting property, for all but finitely many $t\in \mathcal{T}$ and all $\alpha \in
 \mathcal{A}$,  we need to ensure the existence of a collection $\mathcal{C}_{t,\alpha}$ of balls $\mathfrak{B}$ centred in $\overline{B}$ and an appropriate
 sequence $\{k_t\}_{t\in \mathcal{T}}$ of positive numbers satisfying (\ref{contract1}),
 (\ref{contract2}) and (\ref{contract3}) with $\phi=\phi_{\delta}$ and $\phi^{+}=\phi_{\delta}^{+}$.  With this in mind,  let  $(t,\alpha)\in \mathcal{T}\times \mathcal{A} $ and suppose that  $ I_t (\alpha, \phi_{\delta}(t))=\emptyset$.  Then the collection    $\mathcal{C}_{t,\alpha}=\emptyset$ obviously suffices.   Thus, we  can assume that $ I_t (\alpha, \phi_{\delta}(t))  \neq \emptyset$ and in view of \eqref{inhomsetmodi}, it follows that

 \begin{equation}\label{slv22}
 I_t (\alpha, \phi_{\delta}^{+} (t))\cap B   \ \subseteq \ \left \{\bx\in B: |\hat{\theta}(\bx)+ a_0+(\bx, \tilde{\bx}A)\ba|<
 \frac{1}{2^{\left(n-\frac{\delta+\gamma}{2}\right)t}} \right \} \, . \end{equation}
 \noindent Assume for the moment that $\hat\theta$ is a linear map given by \eqref{hatteta}. Then, by \eqref{vb1} and \eqref{diocond+}, we have that $\hat{\theta}(\bx)+ a_0+(\bx, \tilde{\bx}A)\ba$ is a linear combination of $x_1,\dots,x_\xd$ with at least one of the coefficient being $\gg 2^{t(-n+\gamma')}$ in absolute value, where $0<\gamma'<n-\omega(A;\bm\theta)$. Hence,
 $$
 \sup_{\bx\in B}|\hat{\theta}(\bx)+ a_0+(\bx, \tilde{\bx}A)\ba|\gg 2^{t(-n+\gamma')}\,,
 $$
 where the implied constant will not depend on $t$. Therefore, in view of \eqref{slv22}, choosing $\gamma$ within \eqref{eqn:gamma} so that we additionally meet the inequalities
 \begin{equation}\label{eqn:gamma'}
0<\gamma<\gamma'
\end{equation}
ensures that
 \begin{equation}\label{propsubset}
 I_t (\alpha, \phi_{\delta}^{+} (t))\cap B
 \subsetneq B   \qquad \forall  \ \ t\ge t_0 \,,\end{equation}
where $t_0 \in \mathbb{N}$ is a sufficiently large constant.

Now consider the case $\hat\theta$ is not a linear function. Then,
\begin{equation}\label{vb2}
\hat{\theta}(\bx) + a_0 + (\bx,
\tilde{\bx}A)\ba=\hat\theta(\bx)+\tilde{\bx}(A \ba' + \ba'')
\end{equation}
where $\ba'=(a_{\xd+1},\dots,a_n)^t$ and $\ba''=(a_0,\dots,a_\xd)^t$. Thus, \eqref{vb2} is a linear combination of the functions $1$, $x_1$, \dots, $x_\xd$, $\hat\theta(\bx)$ which are linearly independent over $\R$. Therefore, \eqref{vb2}  is not identically zero. Furthermore, the vector of the coefficients of this linear combination is obviously of norm at least $1$. Hence
\begin{align}
\nonumber\inf_{(\ba,a_0)\in\bbR^{n+1}}\sup_{\bx\in B}|\hat{\theta}(\bx) & + a_0 + (\bx,\tilde{\bx}A)\ba|\\[1ex]
&\ge \inf_{\|\bm\eta\|=1}\sup_{\bx\in B}|\eta_0+\eta_1x_1+\dots+\eta_\xd x_\xd+\eta_{\xd+1}\hat{\theta}(\bx)|>0\,,\label{vb3}
\end{align}
where $\bm\eta=(\eta_0,\dots,\eta_{\xd+1})$ and the latter quantity is strictly positive since we take the infimum of a positive continuous function over a compact set (the unit sphere). By \eqref{slv22} and \eqref{vb3}, we once again ensure that \eqref{propsubset} holds for a sufficiently large choice of $t_0$.

 By \eqref{slv1} and the fact that
 $I_t (\alpha, \phi_{\delta}^{+} (t))$ is open, for any $\bx \in I_t (\alpha, \phi_{\delta}(t))\cap \overline{B}$,
 there is  a ball $\mathfrak{B}'(\bx)$
 centred at $\bx$ such that
 \begin{equation} \label{slv3}
 \mathfrak{B}'(\bx)\subseteq I_t (\alpha, \phi_{\delta}^{+} (t))\,.
 \end{equation}

 \noindent On combining  (\ref{propsubset}),  (\ref{slv3})  and the fact that $B$ is bounded, we find that there exists a scaling factor $\tau\ge1$ such that
the ball $\mathfrak{B}(\bx) :=\tau \mathfrak{B}'(\bx)$ satisfies
\begin{equation}\label{contain} \mathfrak{B}(\bx)  \cap   \overline{B}  \ \subseteq  \  I_t (\alpha, \phi_{\delta}^{+} (t))  \cap   \overline{B}   \   \supsetneq \
5\mathfrak{B}(\bx)  \cap \overline{B}   \,\end{equation}
and
\begin{equation}\label{containsv}
5\mathfrak{B}(\bx)   \ \subset  \  11 \, B
\, .\end{equation}
For $t\geq t_0$ and $\alpha\in\mathcal{A}$, we now let
\[\mathcal{C}_{t,\alpha}:=\Big\{\mathfrak{B}(\bx): \bx\in I_t(\alpha,\phi_{\delta}(t))  \cap  \overline{B}  \  \Big\}\,.
\]
Then by construction and the l.h.s. of \eqref{contain}, any such collection of balls  automatically  satisfies conditions (\ref{contract1}) and (\ref{contract2}) with $\phi=\phi_{\delta}$ and $\phi^{+}=\phi_{\delta}^{+}$. Regarding condition (\ref{contract3}), we proceed as follows.\\

Let $\mathfrak{B} \in \mathcal{C}_{t,\alpha}$.  By
(\ref{inhomsetmodi}) and the r.h.s. of (\ref{contain}), we have  that

\begin{equation}\label{obs2}
	\displaystyle \sup_{\bx\in 5\mathfrak{B}}\mathbf{F}_{t,\alpha}(\bx)\geq \sup_{\bx\in 5\mathfrak{B}\cap \overline{B}}
	\mathbf{F}_{t,\alpha}(\bx)\geq \sqrt{n\xd L}\ 2^{\frac{\delta+\gamma}{2} t}\ 2^{t/2}\,.
\end{equation}
On the other hand,
\begin{equation}\label{obs1}
	\displaystyle \sup_{\bx \in 5\mathfrak{B}\cap I_t(\alpha, \phi_\delta (t))} \mathbf{F}_{t,\alpha}(\bx)\leq \sqrt{n\xd L}
	\ 2^{\delta t}\ 2^{t/2}   \, .
\end{equation}
On combining  (\ref{obs2}) and (\ref{obs1}), it follows that
\[\displaystyle \sup_{\bx \in 5\mathfrak{B}\cap I_t(\alpha, \phi_\delta (t))} \mathbf{F}_{t,\alpha}(\bx)
 \ \leq \  2^{\delta t}\ \frac{1}{2^{\frac{\delta+\gamma}{2} t}} \  \sup_{\bx\in 5\mathfrak{B}}\mathbf{F}_{t,\alpha}(\bx)
\ = \ \frac{1}{2^{\frac{\gamma-\delta}{2} t}}\  \sup_{\bx\in 5\mathfrak{B}}\mathbf{F}_{t,\alpha}(\bx)\,. \]
 This together with \eqref{containsv} and (\ref{goodeq}),  implies that for any
 $t\geq t_0$ and $\alpha\in\mathcal{A}$

\begin{eqnarray}\label{lasteq} \varrho(5\mathfrak{B}\cap I_t(\alpha, \phi_\delta (t))) & \leq   &  |5\mathfrak{B}\cap I_t
(\alpha, \phi_\delta (t))| \nonumber \\[3ex]
& \leq   &  \left|\left \{\bx \in 5\mathfrak{B}:\mathbf{F}_{t,\alpha}(\bx)\leq \frac{1}{2^{\frac{\gamma-\delta}{2} t}} \
\sup_{\bx\in 5\mathfrak{B}}\mathbf{F}_{t,\alpha}(\bx) \right\}\right|  \nonumber \\[3ex]
& \leq   &  \tilde{C} \  \frac{1}{2^{\frac{\gamma-\delta}{2} \alpha_0t}} \ \ |5\mathfrak{B}|\,,\end{eqnarray} \\
where $\tilde{C} > 0 $ is some  constant depending on $B$. On using \eqref{containsv} and  the fact that $\mathfrak{B}$ is centred in $\overline{B}$, we have that $|5\mathfrak{B}|\leq c_\xd \varrho(5\mathfrak{B})$ for some
constant $c_\xd$ depending on $\xd $ only.  Hence   (\ref{lasteq}) implies that for all but finitely many $t\in \mathcal{T}$
\begin{eqnarray*}\label{lasteqsv} \varrho(5\mathfrak{B}\cap I_t(\alpha, \phi_\delta (t))) & \leq   &  c_\xd  \tilde{C} \  2^{-\frac{\gamma-\delta}{2} \alpha_0t}   \  \ \varrho(5\mathfrak{B})  \,.\end{eqnarray*}
This verifies condition (\ref{contract3}) with $\phi=\phi_{\delta}$ and
$$
k_t := c_\xd  \tilde{C} \  2^{-\frac{\gamma-\delta}{2} \alpha_0t}  \, .  $$ Furthermore, it is easily seen that  $\sum_{t\in \mathcal{T}}k_t<\infty$ and thus all the conditions of the contracting property are satisfied for the collection $\mathcal{C}_{t,\alpha}$  as defined above.

\section{The divergence theory: proof of Theorem \ref{main divergence} }\label{section: divergence}
The proof of  Theorem \ref{main divergence} makes use  of the  following statement, which is a special and
simplified version of Theorem 3 appearing in  \cite{BaBeVe}.
\begin{theorem}\label {BaBeVe simple}
 Let $\mathcal{M}:= \{ \mathbf{f}(\bx)  : \bx \in U  \}  \subseteq \bbR^n$   be a manifold of dimension $\xd $  parameterized by a smooth map $ \mathbf{f}:U\rightarrow \R^n$  defined on a ball $U\subseteq \bbR^\xd$.  Suppose
 there exists an  absolute constant $C_0 \geq 1$ such that for any ball $B$ with $2B\subseteq U$ and any $\kappa\in (0,1)$, we have that
 \begin{equation}\label{nice eqn}
  \left |\left \{\mathbf{x}\in B: \exists \,  (\mathbf{a}, a_0)\in \mathbb{Z}^n\setminus \{\mathbf{0}\} \times \mathbb{Z} \text{ s.t.}
  \left|
\begin{array}{ll}
|a_0 + \mathbf{f}(\mathbf{x})\cdot\ba| < \displaystyle
\frac{\kappa}{Q^n}\\
    \|\ba\|\leq Q
\end{array}\right.\right\}\right| \ \leq \  C_0  \ \kappa \  |B|
 \end{equation} for all sufficiently large $Q$. Let $\psi$ be an  approximation function and $\theta:\bbR^n \rightarrow \bbR$  be a function such that $\theta|_{\mathcal{M}}\in C^2$.  Let $\xs>\xd-1$ and suppose that
 \begin{equation}\label{eq:divcond+}
 \displaystyle \sum_{\mathbf{a}\in \mathbb{Z}^n\setminus \{\mathbf{0}\}}
 \|\ba\|\left(\frac{\psi (\|\mathbf{a}\|^n)}{\|\ba\|}\right)^{\xs+1-\xd}=\infty\,.
\end{equation}
 Then
\begin{equation}\label{eqn: hau}
 \mathcal{H}^\xs \big(\mathcal{W}^{\theta}_{n}(\psi) \cap \mathcal{M} \big)= \mathcal{H}^\xs(\mathcal{M})\,.
 \end{equation}
\end{theorem}

Note that, by the monotonicity of $\psi$, \eqref{eq:divcond+} is equivalent to \eqref{eq:divcond}. Hence, armed with Theorem \ref{BaBeVe simple}, the proof of  Theorem \ref{main divergence} reduces to establishing \eqref{nice eqn} with   $U\subset \bbR^\xd$ being an open subset, $\mathcal{M}=\mathscr{H}$  and  $\mathbf{f}$ given by \eqref{defaffine}. With this in mind, for any ball  $B$ such that $11B\subseteq U$, any $\kappa\in (0,1)$ and  $Q>1$,  let
\begin{equation}\label{big gra}
 \mathcal{L}^1(B,\kappa, Q):=\left \{\mathbf{x}\in B: \exists (\mathbf{a}, a_0)\in \mathbb{Z}^n\setminus \{\mathbf{0}\} \times \mathbb{Z} \text{ s.t.}
  \left|
\begin{array}{rcl}
|a_0 + (\bx, \tilde{\bx}A)\cdot\ba| <
\frac{\kappa}{Q^n}\\[2ex]
\|\nabla (\bx, \tilde{\bx}A)\cdot \ba\|\geq \frac{\sqrt{n\xd
\|\ba\|}}{2r}\\[2ex] \|\ba\|\leq Q
\end{array}\right.\right\}
\end{equation} and
\begin{equation}\label{small gra}
 \mathcal{L}^2(B,\kappa, Q):=\left \{\mathbf{x}\in B: \exists (\mathbf{a}, a_0)\in \mathbb{Z}^n\setminus \{\mathbf{0}\} \times \mathbb{Z} \text{ s.t.}
  \left|
\begin{array}{rcl}
|a_0 + (\bx, \tilde{\bx}A)\cdot\ba| <
\frac{\kappa}{Q^n}\\[2ex]
\|\nabla (\bx, \tilde{\bx}A)\cdot \ba\|< \frac{\sqrt{n\xd
\|\ba\|}}{2r}\\[2ex] \|\mathbf{a}\|\leq Q
\end{array}\right.\right\}\, .
\end{equation}
We note that the set appearing in (\ref{nice eqn}) is contained in the union of the ``large derivative'' set $\mathcal{L}^1(B,\kappa, Q)$  and  the ``small derivative'' set $\mathcal{L}^2(B,\kappa, Q)$.  Thus
 \begin{equation} \label{lab}
 {\rm l.h.s. \ of \ \eqref{nice eqn}}  \ \leq \ |\mathcal{L}^1(B,\kappa,Q)| \ + \ |\mathcal{L}^2(B,\kappa,Q)|  \, .
\end{equation}
As in the proof of Theorem \ref{main convergence}, estimating the measure of the large derivative set  is relatively easy and makes use of Proposition \ref{prop:large}.  To begin with, observe that
$$\mathcal{L}^1(B, \kappa, Q) \ = \ \displaystyle\bigcup_{\ba \in \mathbb{Z}^n, \,0<\|\ba\|\leq Q}\mathcal{L}^1(\ba,B, \kappa, Q) \, ,  $$ where
for any $\ba\in \mathbb{Z}^n
\setminus \{\mathbf{0}\}$, $$\mathcal{L}^1(\ba,B,\kappa,Q):= \left \{\mathbf{x}\in B:
  \left|
\begin{array}{rcl}
|a_0 + (\bx, \tilde{\bx}A)\cdot\ba| <
\frac{\kappa}{Q^n}\\[2ex]
\|\nabla (\bx, \tilde{\bx}A)\cdot \ba\|\geq \frac{\sqrt{n\xd
\|\ba\|}}{2r}
\end{array}\right.\right\}.$$   Now fix $\ba \in \bbZ^n \setminus \{\mathbf{0}\}$ and with reference to Proposition \ref{prop:large}, let
$F(\bx)=(\mathbf{x},\tilde{\mathbf{x}}A)\cdot \ba$
for $\bx \in 2B$ and $\delta'= \kappa\, Q^{-n}$. By definition,
$$
M \ \geq \
 \frac{1}{4r^2}  \,
$$
where $M$ is given \eqref{defM} and so it follows from Proposition \ref{prop:large} that
$$ | \mathcal{L}^1(\ba,B,\kappa,Q) | \
\leq \  K_\xd \ \displaystyle \frac{\kappa}{Q^n}  \ |B|    \, . $$
In turn, this implies that
\begin{equation}\label{bg estimate}
 |\mathcal{L}^1(B,\kappa,Q)|   \ \leq  \  K_\xd \displaystyle \ \frac{\kappa}{Q^n} (2Q+1)^n  \ |B| \ \leq  \   3^n K_\xd  \ \kappa \  |B|\,.
\end{equation}  \\

We now turn our attention to estimating the measure of the small derivative set  $\mathcal{L}^2(B,\kappa,Q)$. Let $t$ be the unique integer satisfying $2^t\le Q<2^{t+1}$ and $\delta$ satisfy $0<\delta<\gamma$, where $\gamma$ is as defined earlier. For $t$ sufficiently large we obviously have that
\begin{equation}\label{slv3331}
\mathcal{L}^2(B,\kappa,Q)  \  \subset \   \mathcal{A}_t\,,
\end{equation}
where $\mathcal{A}_t$ is defined at the beginning of \S\ref{homo}. Then, as a result of \eqref{subset} and \eqref{vb4}, we have that
\begin{equation}\label{sm estimate}
 |\mathcal{L}^2(B,\kappa,Q)| \ \leq  \  \kappa\  |B|
\end{equation}
provided that $t$ is sufficiently large.

\noindent The desired estimate \eqref{nice eqn}  now follows from \eqref{lab}, (\ref{bg estimate}) and (\ref{sm estimate}) with
$$
C_0 :=  3^n K_\xd+1\, .
$$
This thereby completes the proof of Theorem \ref{main divergence}.

\vspace{7ex}

\noindent{\bf Acknowledgements.}  SV would like to thank the  Tata Institute of Fundamental Research (Mumbai) for its hospitality during a productive and most enjoyable visit during March 2017. On a more personal note, many thanks to Anish and Dhru for sharing with me  the wonderful culinary delights of Mumbai and the wonderful wedding party.  Finally a great big thanks to the lovely  Bridget and the feisty young women Iona and Ayesha for a wonderful 2017; especially the weeks during the summer in France!

\bibliographystyle{amsplain}

\end{document}